\documentclass[journal]{IEEEtran}
\usepackage{lipsum}
\usepackage{amsfonts}

\usepackage{amsmath}
\usepackage[colorlinks,
            linkcolor=cyan,
            anchorcolor=green,
            citecolor=cyan
            ]{hyperref}
\usepackage{cleveref} 
\usepackage{algorithm}
\usepackage[algo2e]{algorithm2e}
\usepackage{float}
\usepackage{graphicx}
\usepackage{epstopdf}
\usepackage{multirow}
\usepackage{caption,subcaption}
 
\def\F{{\mathcal F}}

\def\A{{\mathcal A}}
\def\P{{\Pi}}

\def\R{{\mathbb R}}
\def\T{{ \Theta}}

\def\bfh{{\rm h }}
\def\bfl{{\rm \lambda }}
\def\var{{\boldsymbol  \varepsilon }}
\def\bfa{{\boldsymbol  \phi }}

\def\beps{{\boldsymbol \varepsilon }}
\def\bxi{{\boldsymbol  \xi }}

\def\bfb{{\bf b}}
\def\bfc{{\bf c}}

\def\bfv{{\bf v}}
\def\bfx{{\bf x}}
\def\bfy{{\bf y}}
\def\bfz{{\bf z}}
\def\bfw{{\bf w}}
\def\bfu{{\bf u}}
\def\bfe{{\bf 1}}
\def\bfh{{\bf h}}
\def\supp{{\rm supp}}
\def\Diag{{\rm Diag}}
\def\gpa{{\texttt{GPSP}}}
\def\tr{{\rm true}}

\newcommand{\sgn}{{\mathrm{sgn}}}

\newcommand{\ts}{$\tau$-stationary point}

\newtheorem{theorem}{Theorem}[section]

\newtheorem{lemma}{Lemma}[section]

\newtheorem{corollary}{Corollary}[section]

\newtheorem{remark}{Remark}[section]

\newtheorem{example}{Example}[section]
\graphicspath{ {./images/} }
\ifCLASSINFOpdf
\else
\fi
\begin{document}
%
\title{Computing One-bit Compressive Sensing via Double-Sparsity Constrained Optimization}
%
%
%

\author{Shenglong Zhou,
        Ziyan Luo,
        Naihua Xiu, {and Geoffrey Ye Li, \textit{Fellow, IEEE}}
\thanks{This work is supported in part by the National
Natural Science Foundation of China (12131004) and Beijing Natural Science Foundation (Z190002).}
\thanks{S.L. Zhou and G.Y. Li are with the ITP Lab,  Department of Electrical and Electronic Engineering, Imperial College London, London SW72AZ, United Kingdom,  e-mail: \{shenglong.zhou, geoffrey.li\}@imperial.ac.uk.}
\thanks{Z.Y. Luo and N.H. Xiu are with the Department of Applied Mathematics, Beijing Jiaotong University, Beijing 100044,  People's Republic of China, e-mail: \{zyluo, nhxiu\}@bjtu.edu.cn}
}

%
%

\markboth{}
{Shell \MakeLowercase{\textit{et al.}}: Bare Demo of IEEEtran.cls for IEEE Journals}
%



\maketitle

\begin{abstract}
One-bit compressive sensing gains its popularity in signal processing and communications due to its low storage costs and low hardware complexity.   However, it has been a challenging task to recover the signal only by exploiting the one-bit (the sign) information.  In this paper, we appropriately formulate the one-bit compressive sensing into a double-sparsity constrained optimization problem. The first-order optimality conditions for this nonconvex and discontinuous problem are established via the newly introduced $\tau$-stationarity, {based on which, a gradient projection subspace pursuit (\texttt{GPSP}) algorithm is developed. It is proven that \texttt{GPSP} can converge globally and terminate within finite steps. Numerical experiments have demonstrated its excellent performance in terms of a high order of accuracy with a fast computational speed.}
\end{abstract}


%
\IEEEpeerreviewmaketitle

\section{Introduction}
%
%
%
%
\IEEEPARstart{C}{ompressive} sensing (CS) has seen evolutionary advances in theory and algorithms in the past few decades since introduced in the ground-breaking papers \cite{candes2005decoding,candes2006robust,donoho2006compressed}.
It aims to reconstruct a sparse signal $\bfx$ from an underdetermined linear system $\Phi \bfx=\bfb$, where $\Phi\in\R^{m\times n}$ is the measurement matrix and $\bfb\in\R^{m}$ is the measurement observation. To reduce storage costs and hardware complexity, in \cite{BB08}, only the sign information of measurements are preserved, that is
\begin{eqnarray}\label{exact-recovery-1bcs}
\bfc=\sgn( \Phi\bfx).
\end{eqnarray}
Here,  $\sgn(t)$  returns one if $t$ is positive and negative one otherwise, and thus $c_i\in\{1,-1\}, i\in [m]:=\{1,2,\cdots,m\}$ is the one-bit measurement.  This gives rise to the one-bit CS. It was then extensively applied into applications including communications \cite{haboba2011architecture, movahed2014iterative,  tang2017low, zhou2017sparse},
  wireless sensor network  \cite{meng2009sparse, feng2009multiple, xiong20141, chen2015amplitude}, cognitive radio \cite{lee2012spectrum, fu2014sub}, {imaging science \cite{bourquard2012binary, dong2015map, marcos2016compressed} and to name a few. We refer to a couple of nice surveys \cite{li2018survey,qin2018sparse} for more applications.}
\subsection{Related work}
{\bf a) Noiseless recovery.}  
The task of one-bit CS constructs the sparse signal from the one-bit measurements. The ideal optimization model for recovery problem \eqref{exact-recovery-1bcs} is the following  {$l_0$-norm} minimization,
\begin{eqnarray}\label{1-bit-cs-l0}
\underset{\bfx\in \R^{n}}{\min}~\|\bfx\|_0,~~
\mbox{s.t.} ~ \bfc=\sgn( \Phi\bfx),~ {\|\bfx\|=1},
\end{eqnarray}
where $\|\bfx\|_0$ counts the number of  non-zero entries of $\bfx$ and $\|\bfx\|$ is its $l_2$-norm (Euclidean norm). We note that $\|\bfx\|_0$ is sometimes called {$l_0$-}norm, but it is actually not a norm in a strict mathematical sense.  
 An impressive body of work has developed numerical algorithms for solving problem  \eqref{1-bit-cs-l0}, but most of them focused on its approximations due to the NP-hardness. The earliest work can be traced back to \cite{BB08} where  model \eqref{1-bit-cs-l0} was relaxed by 
\begin{eqnarray}\label{1-bit-cs}
\underset{\bfx\in \R^{n}}{\min}~\|\bfx\|_1,~~
\mbox{s.t.} ~ A \bfx \geq {0},~ \|\bfx\|=1.
\end{eqnarray}
Here, $A:=\Diag(\bfc) \Phi,~ \Diag(\bfc) $ represents the diagonal matrix with diagonal entries from $\bfc$,  and $\|\cdot\|_1$ is the $l_1$-norm. {Other relevant work on the case of exact recovery \eqref{exact-recovery-1bcs} includes the support recovery algorithm \cite{gopi2013one}, the binary iterative re-weighted method \cite{wanghuang2015}, the superset technique approximation \cite{flodin2019superset}, the
fixed-point continuation algorithm \cite{  xiao2019one}, and model-based deep learning \cite{ khobahi2020model}.}

{\bf b) Noisy recovery.} In reality,  measurement $\Phi \bfx$ is frequently contaminated by noise $\var$ before the quantization, i.e., 
\begin{eqnarray}\label{inexact-recovery-1bcs}\bfc=\sgn( \Phi\bfx + {\var}).
\end{eqnarray}
For this scenario,  a popular approach to reconstruct the signal benefits from  the following optimization
\begin{eqnarray}\label{1-bit-cs-loss}
\underset{\bfx\in \R^{n}}{\min}~ \|\bfx\|_1 + \lambda \varphi(\bfx),~~\mbox{s.t.} ~  \|\bfx\|=1,
\end{eqnarray}
where $\lambda>0$ and $\varphi:\R^n\rightarrow\R$ is a loss function. In \cite{BB08}, they adopted  the one-sided $l_2$ function   $\varphi(\bfx):=\|(-A\bfx  )_+\|^2$ with ${\bfy}_+:=(\max\{y_1,{0}\}, \cdots, \{y_m,{0}\})^\top $ and employed  a renormalized fixed point iteration algorithm. Since the targeted problem is a nonconvex optimization, the convergence result has not been provided. The same problem was also addressed by a restricted step shrinkage algorithm \cite{laska2011trust}, where the generated sequence was proved to converge to a stationary point of the penalty problem if some slightly strong assumptions on the sequence were satisfied.

 Following the work in \cite{BB08},  Boufounos modified compressed sampling matching pursuit (CoSaMP) \cite{needell2009cosamp}, one of the most popular greedy methods in CS, to derive the matching sign pursuit method \cite{boufounos2009greedy}.  It turned out to address the sparsity constrained model,
\begin{eqnarray}\label{1-bit-cs-s}
\underset{\bfx\in \R^{n}}{\min}~ \varphi(\bfx),~~\mbox{s.t.} ~  \|\bfx\|_0=s,~\|\bfx\|=1,
\end{eqnarray}
where $s\ll n$  is a given sparsity level and $\varphi$ is the one-sided $l_2$ function.  Based on the framework of the  famous iterative hard thresholding   algorithm, the modified version binary iterative hard thresholding (BIHT)  was then developed in \cite{jacques2013robust} to solve problem  \eqref{1-bit-cs-s}. Apart from the one-sided $l_2$ function, BITH was also able to process the one-sided $l_1$ function, namely, $\varphi(\bfx)=\|(-A\bfx)_+\|_1$.  It was claimed that with a high probability, the distance between a reconstructed signal by BIHT and the original one can be bounded by a prefixed accuracy if the former quantizes to the same quantization point as the latter. As a consequence, the method enjoys a local convergence property.  
{Recently, the normalized BIHT was investigated in \cite{friedlander2020nbiht}, and  achieves an approximation error
rate optimal up to logarithmic factors with a high probability. Other work on noiseless recovery \eqref{inexact-recovery-1bcs} consists of the convex relaxation \cite{plan2012robust},   the passive algorithm \cite{zhang2014efficient}, the soft
consistency reconstructions  \cite{cai2014soft},  the pinball loss iterative hard thresholding \cite{huang2018pinball}, the  sparse consistent coding algorithm \cite{rencker2019sparse}, nonconvex penalty methods \cite{ huang2018nonconvex}, and  feature selection methods \cite{xu2020feature}.}

{\bf c) Noisy recovery with sign flips.} {When some signs are flipped after quantization, the scenario becomes  
\begin{eqnarray}\label{noise-1bcs}
\bfc={\rm Diag}(\bfh)\sgn( \Phi\bfx + \beps),
\end{eqnarray}
where $h_i\in\{1,-1\}$ satisfying $\|\bfh-{\bf1}\|_0\leq k$ and $\beps\in\R^m$ is the noise before quantization. Here, $k$ is a pre-estimated integer and can be deemed as the upper bound of the number of sign flips. Since there is an impressive body of work developing algorithms to address the above problem, see  \cite{fu2014robust, li2014robust, shi2016methods}, we only review a small portion, which is enough to motivate our work in this paper. }

In \cite{huang2018robust}, authors took advantage of the $l_1$-regularized least squares and built the following model
\begin{eqnarray}\label{1-bit-cs-least}
\underset{\bfx\in \R^{n}}{\min}~  \|\bfx\|_1+\lambda\|\bfc-\Phi\bfx\|^2.
\end{eqnarray}
{We note that this optimization problem does not model the loss of amplitude information in $\Phi\bfx$.} However, as stated in \cite{huang2018robust}, with a high probability, the distance between the solution to the model (up to a constant) and a sparse solution  can be bounded by a prefixed accuracy if the sample size $m$ is greater than a threshold. Then a primal dual active set algorithm was proposed  to solve the above model and proved to converge within one step under two assumptions: the submatrix  of $\Phi$ indexed on the nonzero components of the sparse solution is full row rank and the initial point is sufficiently close to the sparse solution. Therefore, the generated sequence again has a local convergence property. Very recently, the authors in \cite{fan2021robust} replaced the $l_1$-norm in \eqref{1-bit-cs-least} by $\|\bfx\|_p^p:=\sum |x_i|^p (0<p<1)$ to design a weighted primal dual active set algorithm. 

To exploit the information of sign flips (i.e., upper bound $k$),  authors in \cite{yan2012robust} integrated a sparse variable $\bfw$  to problem  \eqref{1-bit-cs-s}. The nonzero components in  $\bfw$  represent the measurements that have sign
flips. The resulting optimization problem is
\begin{eqnarray}\label{aop}
\begin{array}{lcl}
&\underset{\bfx\in \R^{n},\bfw\in\R^m}{\min}& \|(-\Diag(\bfc) (\Phi\bfx+{\bfw}))_+\|_p^p,\\
&\mbox{s.t.} &  \|\bfx\|_0\leq s,~\|\bfx\|=1,~\|\bfw\|_0\leq k,\end{array}
\end{eqnarray}
where $p=1$ or $2$. To tackle the above problem,  an alternating minimization method (adaptive outliers pursuit, AOP) was cast:  solving one variable while fixing the other. However, AOP has been tested to heavily rely on the choice of $k$ and the convergence result remains to be seen. Other work relating to \eqref{aop} includes the noise-adaptive renormalized fixed point iteration approach \cite{movahed2012robust} and  the noise-adaptive restricted step shrinkage \cite{movahed2014recovering}. 

When upper bound $k$ of the number of sign flips is unavailable, a remedy pursues a solution with sign flips as few as possible, which can be fulfilled by the following one-sided $l_0$ function minimization \cite{dai2016noisy},
\begin{eqnarray}\label{0-1-regular}
\underset{\bfx\in \R^{n}}{\min}~ \|(\epsilon {\bf1} -A\bfx)_+\|_0 + \eta \|\bfx\|^2,~~\mbox{s.t.} ~  \|\bfx\|_0\leq s,
\end{eqnarray}
where $\eta$ and $\epsilon$ are given positive parameters. Here, $\epsilon$ is used to majorize the objective function. The first term in the objective function  arises from maximizing a posterior estimation from the perspective of statistics. It returns the number of positive components of $(\epsilon {\bf1}-A\bfx)$ and can be regarded as the number of the sign flips when $\epsilon$ is quite small. Instead of solving the one-sided $l_0$ model directly,  a fixed-point algorithm \cite{dai2016noisy} was created for its approximation,
\begin{eqnarray}\label{0-1-regular-approx}
\begin{array}{cl}\underset{\bfx\in \R^{n},\bfw\in\R^m}{\min}& \|(\epsilon {\bf1} -\bfw)_+\|_0 +  \mu\|\bfw-A\bfx\|^2+\eta \|\bfx\|^2,\\
\mbox{s.t.}&  \|\bfx\|_0\leq s,\end{array}
\end{eqnarray}
where $\mu>0$.  It has shown that the generated  sequence converges to a local minimizer of the approximation problem if the spectral norm of matrix $A$ is bounded by some chosen parameters. However, the relationship between the solution obtained by the method and the one to the original problem \eqref{0-1-regular} has not been well explored. 

{To summarize, all these aforementioned methods either had no convergence guarantees (e.g., \cite{BB08,boufounos2009greedy,yan2012robust,movahed2012robust}) or had convergence results that required more or less assumptions on data $\Phi$ or $A$ (e.g., \cite{ huang2018robust,dai2016noisy}). A natural question is whether there is a proper model based on which the proposed algorithm can converge without any assumptions. 
 Moreover, empirical numerical experiments have demonstrated that algorithms that only exploit information of sparsity  did not render desirable accuracies for the case of noisy recovery with sign flips \eqref{noise-1bcs}. By contrast, as shown by numerical experiments in \cite{yan2012robust, movahed2012robust, movahed2014recovering,  dai2016noisy},  the recovery accuracies have been improved greatly when the information of sign flips was taken into consideration.}

\subsection{Our contributions}
{To eliminate the assumptions on data $\Phi$ or  $A$ for better convergence results, we need a loss function with more pleasant features. Moreover, to enhance the recovery accuracy, it is suggested by the work in \cite{yan2012robust, dai2016noisy} that we should take full advantage of the information of both sparsity and the number of sign flips. } Based on these two aspects,  we formulate one-bit CS problem \eqref{noise-1bcs} as 
the following double-sparsity constrained optimization:
 \begin{eqnarray}\label{ell-0-1}
 \begin{array}{lcl}
&\underset{\bfx\in \R^{n},\bfy\in \R^{m}}{\min}&  \| A\bfx+\bfy-\epsilon {\bf1}\|^2+ \eta \|\bfx\|^2,\\
&\mbox{s.t.}&\|  \bfx \|_0\leq s, ~~\|  \bfy_+\|_0\leq k,
\end{array}
\end{eqnarray}
where $\eta>0$ is a penalty parameter, $s\ll n$ and $k\ll m$ are two integers  {representing the prior information on the upper bounds of the signal sparsity and the number of sign flips, respectively. When penalizing the sign flip constraint in our model, it turns to  \eqref{0-1-regular-approx} with  $\bfy=\epsilon {\bf1}-\bfw$. 

{Now we would like to emphasize the merits of our proposed optimization model \eqref{ell-0-1}. For starters, it is able to deal with scenarios \eqref{exact-recovery-1bcs}, \eqref{inexact-recovery-1bcs}, and \eqref{noise-1bcs}. For the first two cases without sign flips, we just set $k=0$ in model \eqref{ell-0-1}. Moreover, it has a strongly convex quadratic  objective function, which is beneficial  to develop a fast second-order algorithm and establish its convergence results   without any additional assumptions, see Theorems \ref{global-convergence}  and \ref{convergence-rate}. Finally, the new model exploits both the sparsity and the upper bound of the number of sign flips, thereby yielding a high order of accuracy, as shown in Figure \ref{fig:ex1-2-cs-all}. It is worth mentioning that the selection of $k$ is very flexible (see Figure \ref{fig:ex1-cs-k}), which reveals that our approach does not rely on $k$ heavily (while AOP in \cite{yan2012robust} does).}

The main contributions in this paper are threefold:

\begin{itemize}
\item[i.] {\bf The new optimization model.} The double-sparsity constrained optimization, problem  \eqref{ell-0-1}, is formulated to handle the one-bit CS. It is well-known that the  two discrete and nonconvex constraints in \eqref{ell-0-1} lead to the NP-hardness in general. Nevertheless, a necessary and sufficient optimality condition as stated in \eqref{KKT-point} for a local minimizer is established, see \Cref{nec-suff-opt-con-KKT}.   Moreover, the necessary or sufficient optimality condition for a global minimizer is further studied through the newly introduced $\tau$-stationary point, see \Cref{nec-suff-opt-con-sta}. {Finally, it turns out that the distance between any local minimizer to \eqref{ell-0-1} (up to a constant) and the true signal can be bounded by an error bound with a high probability under some assumptions, see Theorem \ref{oracle-property}.}

\item[ii.] {\bf The efficient {\tt GPSP} algorithm.} As the established optimality conditions indicating a $\tau$-stationary point is instructive to pursue an optimal solution to \eqref{ell-0-1}, we design a gradient projection method with a subspace pursuit scheme interpolated (dubbed as {\tt GPSP}). {The proposed method is proved to be globally convergent to a $\tau$-stationary point (denoted by $\bfz^*$), which must be a unique local minimizer of \eqref{ell-0-1} without any assumptions. Moreover,  if we assume an additional condition, then $\bfz^*$ can be a global minimizer, as shown in \Cref{global-convergence}. 
 Furthermore, the produced sequence is eventually identical to  $\bfz^*$, namely, {\tt GPSP} can stop within finite steps, see \Cref{convergence-rate}.}

\item[iii.]{\bf High numerical performance.} {\tt GPSP} is demonstrated to be relatively robust to parameters $k$, $\epsilon$, $\eta$ in \eqref{ell-0-1} in the numerical experiments, which indicates that we do not need an exact upper bound $k$ of the sign flips. 
{In addition, {\tt GPSP} is cast based on problem  \eqref{ell-0-1} and turns out to be a second-order method, thereby leading to a considerably high order of accuracy. Finally, the algorithmic design enables it to have a low computational complexity. 
 Therefore,}  it outperforms {all benchmark  solvers} for synthetic data, in terms of time efficiency and recovery accuracy. 
\end{itemize}}

\subsection{Organization}
The remainder of the paper is organized as follows. In \Cref{sec:preliminary}, some necessary mathematical backgrounds are provided, including the notation and the projection onto the feasible set of problem  \eqref{ell-0-1}.  \Cref{sec:opt}  {is devoted to} the optimality conditions of the problem, associated with the $\tau$-stationary points, followed by its relationship to the global minimizers.  {In \Cref{sec:rec-per}, we investigate the reconstruction performance to the recovery accuracy of solutions to  (\ref{ell-0-1}) and the true signal.
In \Cref{sec:Newton}, the gradient projection subspace pursuit ({\tt GPSP}) method is designed,  {and properties of the global convergence and termination within finite steps are established}.}  Numerical experiments are given in \Cref{sec:numerical}, including the involved parameters tuning and comparisons with other  six excellent solvers.   {Concluding remarks are made in \Cref{sec:conclude}.}

 \section{Preliminaries}\label{sec:preliminary}
\noindent We first define some notation employed throughout this paper. To differ from ${\rm sgn(t)}$, the sign function is written as ${\rm sign(t)}$ that returns $0$ if $t=0$ and ${\rm sgn(t)}$ otherwise. Given a subset $T\subseteq[n]:=\{1,2,\cdots,n\}$, its cardinality and complementary set are $|T|$ and $\overline T:= [n]\setminus T$. For a vector $\bfx \in \R^{n}$, the support set, $\supp(\bfx)$, represents the indices of nonzero elements of $\bfx$ and the neighbourhood  with a radius $\delta>0$ is denoted by  $N( \bfx,\delta):=\{\bfw\in\R^n: \|\bfw- \bfx \|< \delta\}$. Let $\|\cdot\|_{\infty}$ be the infinity norm, $\|\bfx\|_{[i]}$ be the $i$th largest (in absolute) element of $\bfx$, and $\|A\|_2$ be the spectral norm of $A$. In addition, $\bfx_{T}$ stands for the sub-vector  contains elements of $\bfx$ indexed on $T$. Similarly, for a matrix $A\in\mathbb R^{m\times n}$, $A_{\Gamma T}$ is the sub-matrix containing rows indexed on $\Gamma$ and columns indexed on $T$, particularly, $A_{:T}=A_{[m]T}$. Moreover, we merge two vectors  $\bfx$ and $\bfy$ by $\bfz:=(\bfx; \bfy):=(\bfx^\top ~\bfy^\top)^\top$.  For a positive  {definite} matrix $H$, the $H$-weighted norm is written $\|\bfz\|_{H}^2=\langle \bfz, H \bfz\rangle$, where $\langle\bfz,\bfz' \rangle:=\sum z_iz_i'$ is the inner product of two vectors.  {Given a scalar $a\in\R$, $\lceil a\rceil$ returns the smallest integer that is no less than $a$.}  {For simplicity, denote}
\begin{eqnarray}
\label{s-k}\begin{array}{lll}
S&:=&\left\{\bfx\in\R^n:~\|\bfx\|_0\leq s\right\},\\
K&:=&\left\{\bfy\in\R^m:~\|\bfy_+\|_0\leq k\right\}.\end{array}
\end{eqnarray}
 {The} feasible region  of \eqref{ell-0-1}  is then denoted by $$\F:=S\times K.$$ 
\subsection{Projections}
For a nonempty and closed set $\Omega\subseteq\R^n$, the projection $\P_\Omega(\bfx)$ of $\bfx\in\R^n$ onto $\Omega$ is given by $$
\P_\Omega(\bfx) = {\rm argmin} ~  \{ \|\bfx-\bfw\|: \bfw\in\Omega  \}.$$
By introducing
\begin{equation}\label{psx-T}
\Sigma(\bfx;s):=\left\{T\subseteq [n]:\begin{array}{l}
|T|=s,\\
 |x_i|\geq |x_j|, {\forall i\in T},j\notin T 
\end{array} \right\},
\end{equation}
 {one can easily verify that}
\begin{equation}\label{psx}
\P_S(\bfx)=   \Big\{\left(\bfx_{{T}};0\right):~T\in \Sigma(\bfx;s) \Big\}.
\end{equation}
To derive the projection of a point $\bfy \in \R^{m}$  onto $K$, denote
\begin{eqnarray}\label{notation-z} 
\Gamma_+&:=&\left\{i\in[m]:~ y_i>0  \right\},\nonumber\\
\Gamma_0&:=&\left\{i\in[m]:~y_i=0 \right\},\\
\Gamma_-&:=&\left\{i\in[m]:~ y_i<0 \right\}. \nonumber
\end{eqnarray}
Note that $\Gamma_+, \Gamma_0$ and $\Gamma_-$ should depend  on $\bfy$. We drop their dependence if no extra explanations are provided for the sake of notational convenience. Based on the above notation, for a point $\bfy \in \R^{m}$ and an integer $k\in [m]$,  we define a set by
\begin{eqnarray}
\label{T-z}
\arraycolsep=0pt\def\arraystretch{1}
\begin{array}{lll}
&&\T(\bfy;k)\\
&:=&\left\{{\Gamma_{k} \cup  \Gamma_-}: 
\begin{array}{l}
{\Gamma_{k}\subseteq  \Gamma_+}, |\Gamma_{k}|=\min\{k, |\Gamma_+|\} \\
y_i\geq y_j\geq0,{\forall i\in \Gamma_{k}},\forall j\in   \Gamma_+\setminus\Gamma_{k} 
\end{array}
\right\},
\end{array}
\end{eqnarray}
where $\Gamma_+$ and $\Gamma_-$ are given by \eqref{notation-z}. One can observe that $\Gamma\in\T(\bfy;k)$ consists of the indices of all negative elements and the first $\min\{k, |\Gamma_+|\}$ largest positive elements of $\bfy$.  These notation  allow us to derive projection $\P_K(\bfy)$ by
\begin{equation}\label{pky}
\P_K(\bfy)=   \Big\{\left(\bfy_{\Gamma};0\right):~\Gamma\in\T(\bfy;k)   \Big\}.
\end{equation}
For example,  $\bfy=(3,2,2,0,-2)^\top$, then we have
\begin{eqnarray*}
\T(\bfy;3)&=& \Big\{\{1,2,3,5\} \Big\},~~\P_K(\bfy)= \{\bfy \},\\
\T(\bfy;2)&=& \Big\{\{1,2,5\},\{1,3,5\} \Big\},\\
\P_K(\bfy)&=& \Big\{(3,2,0,0,-2)^\top, (3,0,2,0,-2)^\top \Big\}.
\end{eqnarray*}
\subsection{Properties of the objective function}
To end this section, we present some properties of the objective function in \eqref{ell-0-1}, which can be written as follows
\begin{eqnarray}\label{ell-0-1-obj}
f(\bfx,\bfy)&:=&   \| A\bfx+\bfy-\epsilon {\bf1}\|^2+ \eta \|\bfx\|^2 \nonumber\\
&= &  \|  \bfz\|^2_H  - 2\epsilon\langle (A^\top {\bf1};{\bf1}),\bfz \rangle +  m\epsilon^2 \\
 &=:&f(\bfz),\nonumber
\end{eqnarray}
where $H$ is given by
  $$H:=\frac{1}{2}\nabla^2f(\bfz)=\left[
  \begin{array}{ll}
  A^\top A+\eta I & A^\top\\
  A & I
  \end{array}
  \right].$$
 {It is easy to verify that $H$ is symmetric positive definite and hence has all eigenvalues positive. Denote the smallest and the largest eigenvalues by $\lambda_{\min}$ and $\lambda_{\max}$, respectively. The quadratic objective function, $f$, is then strongly convex and strongly smooth since}
for any $\bfz$ and $\bfz'$ in $\R^{n+m}$,
\begin{eqnarray}\label{ell-0-1-obj-convex-smooth}
 &&f(\bfz)-f(\bfz')-\langle \nabla f(\bfz'), \bfz-\bfz'\rangle
\\ &=&  {\|  \bfz-\bfz'\|^2_H}~ \in~ \left[{\lambda_{\min} \|  \bfz-\bfz'\|^2},~  {\lambda_{\max} \|  \bfz-\bfz'\|^2} \right]. \nonumber
\end{eqnarray}
\section{Optimality Conditions} \label{sec:opt}
The first-order necessary and sufficient optimality conditions for problem (\ref{ell-0-1}) are established in this section and all proofs are given in Appendix \ref{app:opt}.
\begin{lemma}\label{nec-suff-opt-con-KKT} Consider a point $\bfz^*:=(\bfx^*;\bfy^*)\in\F$ with
$$T^*:=\supp(\bfx^*),~~\Gamma^*:=\supp(\bfy^*).$$
{A point $\bfz^*$ is a local minimizer  of   (\ref{ell-0-1}) if and only if it satisfies
  \begin{eqnarray}\label{KKT-point} \begin{array}{rlr}
  \begin{array}{r} \nabla_{\bfx} f(\bfz^*)= 0,\end{array}  &\rm{if}&\|\bfx^*\|_0<s,\\
  \begin{array}{r} (\nabla_{\bfx} f(\bfz^*))_{T^*} = 0,\end{array} &\rm{if}&\|\bfx^*\|_0=s,\\
 \begin{array}{r}
 (\nabla_{\bfy} f(\bfz^*))_{\Gamma^*} = 0,\\
    (\nabla_{\bfy} f(\bfz^*))_{\overline\Gamma^*} \leq 0,
    \end{array}
  &\rm{and}&\|\bfy^*_+\|_0=k.  
\end{array} 
   \end{eqnarray}}
Furthermore, for any local minimizer $\bfz^*$, there is a $\delta_*>0$ satisfying  the following quadratic growth property
\begin{eqnarray} \label{dw-dw}
f(\bfz)-f( \bfz^*) \geq {\|\bfz - \bfz^*\|^2_{H}},~~\forall~\bfz\in\F\cap N(\bfz^*, \delta_*).
\end{eqnarray}
\end{lemma}
\begin{remark}\label{remark-exlude-zero}
If  $\bfx^*=0$ is the optimal solution to problem (\ref{ell-0-1}), then $\bfy^*=(\bfy^*_{\Gamma^*};0)=(\epsilon{\bf1};0)$ for any $|\Gamma^*|=k$.  From Lemma \ref{nec-suff-opt-con-KKT},  $\nabla_{\bfx} f(\bfz^*)= -2\epsilon A_{\overline \Gamma^*:}^\top {\bf1}=0$. Therefore, if we need to exclude zero solution, we assume $\|A_{ \Upsilon:}^\top {\bf1}\|_\infty>0$ for any $|\Upsilon|=m-k$. This is a very weak assumption.
\end{remark} 

Lemma \ref{nec-suff-opt-con-KKT} shows the optimality conditions of a point being a local minimizer. We further  establish the conditions for a global minimizer. To do that, we introduce a $\tau $-stationary point.  A point $\bfz^*:=(\bfx^*; \bfy^*)$  is called a $\tau $-stationary point of  (\ref{ell-0-1}) with some $\tau >0$ if it satisfies  \begin{eqnarray}\label{eta-point-z}\bfz^* \in  \P_\F\left( \bfz^* - \tau   \nabla f(\bfz^*)\right). \end{eqnarray}
 {An equivalent characterization is presented as follows.}
  \begin{lemma}\label{pro-eta}  A point $\bfz^*$ is a  $\tau $-stationary point of  problem  (\ref{ell-0-1}) with some $\tau >0$ if and only if it satisfies 
   \begin{eqnarray}  \label{eta-point-1}
  \begin{array}{l}
   \|\bfx^* \|_0\leq s,\\
 \tau(\nabla_{\bfx} f(\bfz^*))_i\left\{
\begin{array}{ll}
=0,&i\in T^*,\\
\in [-\|\bfx^* \|_{[s]},\|\bfx^*\|_{[s]}],&i\in \overline T^*,
\end{array}
\right.\\
{\|\bfy^*_+\|_0 = k},\\
 \tau (\nabla_{\bfy} f(\bfz^*))_i\left\{
\begin{array}{ll}
=0,&~~~~~~i\in \Gamma^*,\\
\in [-\|\bfy^*_+\|_{[k]},0 ],&~~~~~~i\in\overline\Gamma^*.
\end{array}
\right.
\end{array}
  \end{eqnarray}
\end{lemma}

 The following theorem reveals the relationships between  $\tau $-stationary points and global minimizers of problem  (\ref{ell-0-1}).
\begin{theorem}\label{nec-suff-opt-con-sta} 
{For problem (\ref{ell-0-1}) and a point $\bfz^*\in\F$, a global minimizer $\bfz^*$ is   a $\tau $-stationary point with $0<\tau \leq  1/(2\lambda_{\max})$, and conversely, a $\tau $-stationary point with $\tau \geq 1/(2\lambda_{\min})$ is also a global minimizer.}
\end{theorem}

{Based on the above theorem, we can check if a local minimizer is a global minimizer explicitly with the help of the $\tau $-stationary point.}
{
\begin{corollary}\label{check-global} 
Let $\bfz^*$ be a local minimizer of   problem (\ref{ell-0-1}).   Then it  is a  $\tau_*$-stationary point   with
\begin{eqnarray}  \label{eta-point-global-tau*}
\tau_*:=  \left\{\begin{array}{lll}
\frac{\|\bfy^*_+\|_{[k]}}{\|\bxi^*\|_{\infty}}, &  \|\bfx^* \|_0< s,\\
 \min\left\{\frac{\|\bfy^*_+\|_{[k]}}{\|\bxi^*\|_{\infty}},\frac{\|\bfx^*\|_{[s]}}{\|A_{\overline\Gamma^*\overline T^*}^\top\bxi^*\|_{\infty}}\right\}, &  \|\bfx^* \|_0 = s,\\ 
\end{array}\right.
  \end{eqnarray}
where $\bxi^*:=2\epsilon\eta(A_{\overline\Gamma^*T^*}A_{\overline\Gamma^*T^*}^\top+\eta I)^{-1}{\bf1} $. Moreover, $\bfz^*$ is also a global minimizer of problem (\ref{ell-0-1}) if $\tau_* \geq 1/(2\lambda_{\min})$. 
\end{corollary}}
{We give a simple example to illustrate the above corollary.
\begin{example}\label{ex-check-global} Consider problem (\ref{ell-0-1}) with $s=1, k=1$ and $A\in\R^{3\times 4}$ given by
\begin{eqnarray}  \label{ex-A}
A:=  \left[\begin{array}{rccc}
 -1& $t$& $t$& 0\\
  1& $t$& 0& 0\\
   1& 0& $t$& $t$\\ 
\end{array}\right],
  \end{eqnarray}
  where $t>0$. Let $c:=\frac{2\epsilon}{2+\eta}$. Now for point $\bfz^*=(\bfx^*;\bfy^*)$ with $\bfx^*= (c,0,0,0 )^\top$ and $
\bfy^*= (c+\epsilon,0,0 )^\top$,
 one can check that $\nabla_{\bfy} f(\bfz^*)=  
   -c\eta (0,
  1,
  1)^\top$ and $
   \nabla_{\bfx} f(\bfz^*)=   -tc\eta  (0,1,1,1)^\top.$
   Therefore, $(\bfx^*, \bfy^*)$ is a local minimizer of problem (\ref{ell-0-1}) since it satisfies \eqref{KKT-point}. Moreover, direct calculations can check that  \begin{eqnarray*}
 \begin{array}{l}\bxi^*=  c \eta{\bf1}~~{\rm and}~~\tau_*= \min\Big\{\frac{4+\eta}{2+\eta},\frac{1}{ \eta t}\Big\}.
  \end{array} \end{eqnarray*} 
  We plot $1/(2\lambda_{\min})$ and $\tau_*$ by  fixing $\eta=1$   or fixing $t=1/5$ in Figure \ref{fig:tau-lam}. One can see that  $\tau_*>1/(2\lambda_{\min})$ when $t\in(0,0.4]$ in  Figure \ref{fig:tau-lam-1} or when $\eta\in[0.1,6]$  in Figure \ref{fig:tau-lam-2}, which by Corollary \ref{check-global} means that $\bfz^*$ is a global minimizer of problem (\ref{ell-0-1}). This example indicates that there are many cases of $A$ and $\eta$ for which the found local minimizer is global.  
\begin{figure}[!th]
\centering
\begin{subfigure}{.24\textwidth}
	\centering
	\includegraphics[width=1\linewidth]{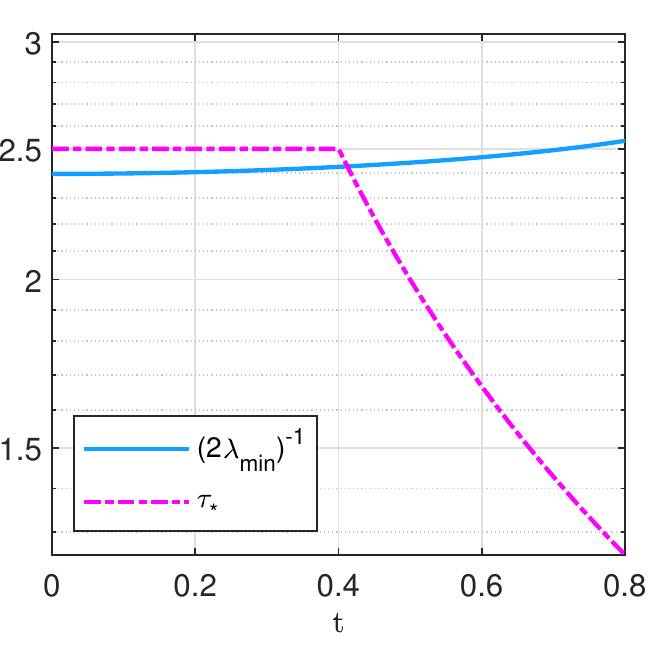}
	\caption{Fixing $\eta=1$.}\label{fig:tau-lam-1}  
\end{subfigure}
\begin{subfigure}{.24\textwidth}
	\centering
	\includegraphics[width=1\linewidth]{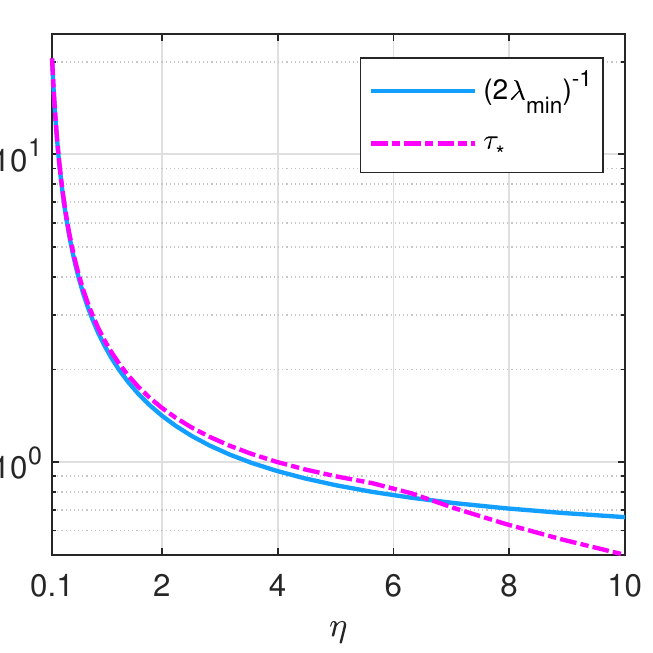}
	\caption{Fixing $t=1/5$.}\label{fig:tau-lam-2}  
\end{subfigure} 
\caption{Values of $1/(2\lambda_{\min})$ and $\tau_*$.}\label{fig:tau-lam}  
\end{figure}
\end{example}
\section{Reconstruction Performance}\label{sec:rec-per}
In this section, we investigate the reconstruction performance of our proposed optimization model (\ref{ell-0-1}). Hereafter, we denote the ground-truth signal by $\bfx^\tr$ with
\begin{eqnarray}  \label{ground-truth}
\|\bfx^\tr\|_0\leq s,\qquad\|\bfx^\tr\|=1.   
  \end{eqnarray}
Moreover, let $\bfz^*=(\bfx^*;\bfy^*)$ with $\bfx^*\neq0$ be any local minimizer of problem (\ref{ell-0-1}) and $$\widehat \bfx : = \bfx^*/\|\bfx^*\|.$$
To exclude zero solution $\bfx^*=0$, as shown in Remark \ref{remark-exlude-zero}, we assume $\|A_{ \Upsilon:}^\top {\bf1}\|_\infty>0$ for any $|\Upsilon|=m-k$ in this section. We now aim at estimating the bound for the gap between $\widehat \bfx $ and $ \bfx^\tr$. To proceed with that, we also need the concept of the binary $\delta$-stable embedding  of a mapping \cite[Definition 1]{jacques2013robust}. Let $\delta\in(0,1)$, a mapping $\A:\R^n\rightarrow\{-1,1\}^m$ is a binary $\delta$-stable embedding (B$\delta$SE) of order $s$ for sparse vectors if
\begin{eqnarray}  \label{B-delta-SE}
d(\bfu,\bfv)-\delta\leq \frac{1}{m}\|\A(\bfu)-\A(\bfv)\|_0\leq d(\bfu,\bfv)+\delta,
  \end{eqnarray}
  for all $\bfu,\bfv$ satisfying $|\supp(\bfu)\cup\supp(\bfv)|\leq s$ and $ \|\bfu\|=\|\bfv\|=1$, where $$d( \bfu,\bfv):=\frac{1}{\pi}{\rm arccos}\langle  \bfu,\bfv\rangle. $$  
The following result states that mapping $\sgn(\cdot)$  is a  B$\delta$SE.  
\begin{lemma}[{\cite[Theorem 3]{jacques2013robust}}] \label{lemma-bse} Let $\delta\in(0,1)$, $\Phi$ be a matrix with entries generated from the independent and identically distributed (i.i.d.) samples of the standard Gaussian distribution $\mathcal{N}(0,1)$  and the number of measurements satisfies 
$$m\geq m(s,\delta,\theta):= \frac{2}{\delta^2}\Big(s \log(n) +2s\log( {35}/{\delta})+\log({2}/{\theta})\Big)$$
for a fixed $\theta\in[0,1]$, then 
\begin{eqnarray*}  
~~~~\mathbb{P}\left\{ |d(\bfu,\bfv)-\|\A(\bfu)-\A(\bfv)\|_0/m |\leq \delta\right\}\geq   1-\theta.
  \end{eqnarray*}
\end{lemma}
\begin{lemma}\label{feasible-recovery}Let $\bfz^*$ be any local minimizer of problem (\ref{ell-0-1}) and $c_A$ represent the minimal eigenvalue of $ A_{\overline \Gamma T}A_{\overline \Gamma T}^\top$ for any $|\Gamma|\geq k$ and $|T|\leq s$. Then for any
\begin{eqnarray}  \label{eta-tuning}
0 \leq \eta \leq  \frac{c_A}{\sqrt{m-k}-1},   
  \end{eqnarray}
point $\widehat{\bfx}$ satisfies 
\begin{eqnarray}  \label{hat-x}
\|\widehat{\bfx}\|=1, ~~ \|\widehat{\bfx}\|_0\leq s, ~~ \|(-A\widehat{\bfx})_+\|_0\leq k.
  \end{eqnarray}
\end{lemma}
\begin{theorem}\label{oracle-property}
Let $\var$ be a noise vector with entries  from i.i.d. samples of  $\mathcal{N}(0,\varrho^2)$ and $\Phi$ be generated as Lemma \ref{lemma-bse} with $m\geq m(2s,\delta,\theta)$ for a fixed $\theta\in[0,1]$ and $\delta\in(0,1)$ satisfying \begin{eqnarray}  \label{gap-x-xtrue-phi}
 \phi:=\frac{k}{m}    + \frac{\varrho}{4}+ \delta\in(0,1/2]. 
  \end{eqnarray}
 Let $\bfz^*$ be any local minimizer of problem (\ref{ell-0-1}) with $\eta$ being chosen as \eqref{eta-tuning}, then 
\begin{eqnarray}  \label{gap-x-xtrue}
\mathbb{P} \{\|\widehat{\bfx}-\bfx^\tr\|  \leq 2\sin(\phi\pi)\} \geq 1-e^{-2m\delta^2}-\theta.
  \end{eqnarray}
\end{theorem}
\begin{remark}\label{remark-para} We have some comments on Theorem \ref{oracle-property}.
\begin{itemize}
\item The result is valid for any local minimizer. If $\phi$ in \eqref{gap-x-xtrue-phi} tends to zero, then we can conclude that  any local minimizer approaches the true signal with a high probability, which indicates there only  exists one local minimizer of problem (\ref{ell-0-1}). This local minimizer is also the unique global minimizer and is the true signal.
\item The reconstruction performance apparently depends on sparsity level $s$, upper bound $k$ as well as number of measurements $m$.  It also provides a hint to set the regularized parameter, $\eta$, by \eqref{eta-tuning}, though the calculation of $c_A$ is quite expensive. Nevertheless, \eqref{eta-tuning} means that $\eta$ should be chosen smaller than a threshold in practical computation. In other words, it should not be set too large, which can be verified by our numerical experiments, see Figure \ref{fig:ex1-cs-eta}.
Moreover, the establishments of Theorem \ref{oracle-property} and Lemma \ref{feasible-recovery}  do not impose assumptions on parameter $\epsilon$. Therefore, the reconstruction performance might be quite robust to   $\epsilon$. This is well testified by our numerical simulation, see Figure \ref{fig:ex1-cs-eps}. 
\end{itemize}
\end{remark}
}
\section{Gradient Projection Subspace Pursuit}\label{sec:Newton}
A gradient projection method with a subspace pursuit strategy  {is proposed to handle problem  \eqref{ell-0-1} by seeking a \ts.}   For notational simplicity, hereafter,  for a parameter $\tau\in(0,1]$, we always let
   \begin{eqnarray}
\label{z-l-p-f}
\bfz^\ell(\tau)=\left[\bfx^\ell(\tau);~\bfy^\ell(\tau)\right]  \in \P_{\F}(\bfz^\ell -\tau \nabla  f(\bfz^\ell)),
  \end{eqnarray}  {for the $\ell$th iteration $\bfz^\ell:=(\bfx^\ell; \bfy^\ell)$.}
 {Analogous to the $\Gamma$-related indices defined for $\bfy$ in \eqref{notation-z}, we} also define
  \begin{eqnarray}\label{notation-y-ell}
  \begin{array}{lll}
\widetilde\Gamma_+^\ell &:=& \{i\in[m]: (\bfy^\ell(\tau_\ell))_i>0\},\\
\Gamma_+^\ell&:=&\left\{i\in[m]:~ y^\ell_i>0  \right\},\\
\Gamma_0^\ell&:=&\left\{i\in[m]:~ y^\ell_i=0  \right\},\\
\Gamma_-^\ell&:=&\left\{i\in[m]:~ y^\ell_i<0  \right\}.
\end{array}\end{eqnarray}
We denote the support sets for $\bfx^\ell$ and $\bfx^\ell(\tau_\ell)$ as follows
  \begin{eqnarray}\label{notation-x-ell}\hspace{10mm}
T^\ell:=\supp(\bfx^\ell),~~ \widetilde T^\ell:=\supp(\bfx^\ell(\tau_\ell)).
\end{eqnarray}
 {Given $\bfz^\ell=(\bfx^\ell;\bfy^\ell)\in\F$, define the following subspace}
\begin{eqnarray}\label{subspace-z}
\Omega(\bfz^\ell):=\left\{\bfz=(
  \bfx;
  \bfy):\begin{array}{l}
\bfx_{\overline{T}^\ell } =0,\\
\bfy_{\Gamma_0^\ell} = 0,~\bfy_{\Gamma_-^\ell} \leq 0
\end{array}\right\}.
\end{eqnarray}
It is easy to see that $\Omega(\bfz^\ell)\subseteq \F$.  
Based on these notation, we summarize the framework of the proposed method in  \Cref{Alg-NL01}.
\begin{algorithm}[!th]
\SetAlgoLined
 Initialize $\bfz^0\in\F,\beta\in(0,1)$ and $\rho,\varepsilon, {\tt tol}_0>0$.  Set $\ell :=0$.\\
 \While{${\tt tol}_\ell > \varepsilon$}{
  \underline{\tt Gradient descent:}\\
   Find the smallest integer $\sigma = 0, 1, 2, \cdots$ such that \\
 \parbox{.45\textwidth}{\begin{eqnarray}\label{armijio-descent-property}f(\bfz^\ell(\beta^\sigma))\leq f(\bfz^\ell) - {\rho}\|\bfz^\ell(\beta^\sigma)-\bfz^\ell\|^2.
\end{eqnarray}}\\
  Set $\tau_\ell=\beta^\sigma$,  $ \bfu^{\ell}:=\bfz^\ell(\tau_\ell)$ and $\bfz^{\ell+1}=\bfu^{\ell}$.\\
 \underline{\tt Subspace pursuit:}\\
 \If{{$(T^\ell= \widetilde T^\ell~\text{or}~\| \nabla_\bfx f(\bfu^\ell)\|\leq \varepsilon)$ and $\Gamma_+^\ell= \widetilde \Gamma_+^\ell$}
}{
   \parbox{.42\textwidth}{\begin{eqnarray}\label{Newton-descent-property}
	\bfv^\ell =  {\rm argmin}~\{f\left(\bfz\right):\bfz\in \Omega(\bfz^\ell)\}.
\end{eqnarray}}\\
   If $f(\bfv^\ell)\leq f(\bfu^\ell) - {\rho}\|\bfv^\ell-\bfu^\ell\|^2$ then set $\bfz^{\ell+1}=\bfv^{\ell}.$
}

Compute ${\tt tol}_\ell:=\|\bfu^{\ell}-\bfz^{\ell}\|$ and set $\ell:= \ell+1$.
 }
 Output solution $\overline{\bfx}=\bfx^\ell/\|\bfx^\ell\|.$
 \caption{ \gpa: Gradient projection subspace pursuit \label{Alg-NL01}}
\end{algorithm}

 {Observing that initial point $\bfz^0 \in \F$, and $\Omega(\bfz^\ell)\subseteq \F$, we can see that all iterations are feasible. Particularly, if gap ${\tt tol}_\ell=\|\bfz^{\ell}-\bfu^\ell\|$ vanishes, then $$\bfz^\ell=\bfu^\ell \in\P_\F( \bfz^\ell -\tau_\ell \nabla f(\bfz^\ell)),$$ which indicates that $\bfz^\ell$ is a $\tau$-stationary point 
with $\tau\leq  \tau_\ell$.}
 {Additionally, once conditions $T^\ell= \widetilde T^\ell$ and $\Gamma_+^\ell= \widetilde \Gamma_+^\ell$ are satisfied, we have $\bfu^l\in \Omega(\bfz^\ell)$. The unique minimizer, $\bfv^\ell$, of $f$ over $\Omega(\bfz^\ell)$ implies that }
\begin{eqnarray}\label{opt-cond-z}
 \langle \nabla  f(\bfv^\ell),  \bfu^\ell-\bfv^\ell  \rangle \geq 0.\end{eqnarray}
In virtue of \eqref{ell-0-1-obj-convex-smooth}, we have 
\begin{eqnarray}\label{newton-step-v}
  f(\bfv^{\ell})&\leq&   f(\bfu^{\ell})- \lambda_{\min} \|\bfu^{\ell}-\bfv^\ell\|^2.
\end{eqnarray}
Suppose $0<\rho\leq \lambda_{\min}$. Then candidate $\bfv^{\ell}$ will be taken, namely, $\bfz^{\ell+1}=\bfv^{\ell}$.

\subsection{Computational complexity analysis}
To update $\bfu^{\ell}$, we need to select one point from $\P_{\F}(\bfz^\ell -\tau \nabla  f(\bfz^\ell))$. Namely, three quantities are computed: $\bar\bfz^\ell:=(\bar\bfx^\ell;\bar\bfy^\ell):=\bfz^\ell -\tau \nabla  f(\bfz^\ell)$, $\P_{S}(\bar\bfx^\ell)$ and $\P_{K}(\bar\bfy^\ell)$. For the former, the computational complexity is about $O(mn)$. To select one point from $\P_{S}(\bar\bfx^\ell)$, we only pick the first $s$ largest (in absolute) elements of $\bar\bfx^\ell$. This allows
us to use a MATLAB  built-in function {\tt maxk}  whose computational
complexity is $O(n + s \log s)$.  {Similarly, for $\P_{K}(\bar\bfy^\ell)$, the computational complexity is $O(m + k \log k)$. Thus, updating $\bfu^{\ell}$ takes a computational complexity of order $O(\sigma mn)$}, 
where $\sigma$ is the smallest integer satisfying (\ref{armijio-descent-property}).

To update $\bfv^{\ell}$, we solve a quadratic programming, \begin{eqnarray} \label{QP}
	\bfv^\ell = \underset{(\bfx;\bfy)}{\rm argmin} &&  \| A\bfx+\bfy-\epsilon {\bf1}\|^2+ \eta \|\bfx\|^2\\
	{\rm s.t.}&& \bfx_{\overline{T}^\ell } =0,~\bfy_{\Gamma_0^\ell} = 0,~\bfy_{\Gamma_-^\ell} \leq 0, \nonumber
\end{eqnarray}
for  fixed $T^\ell$, $\Gamma_0^\ell$ and $\Gamma_-^\ell$. Any solvers for solving the quadratic programming can be used to solve  (\ref{QP}) to pursue a solution in good quality. 
 {To further reduce the computation cost, we drop  constraint $\bfy_{ \Gamma_-^\ell} \leq 0$ from (\ref{QP}) and simply solve the equations:} 
\begin{eqnarray*}
\left[
\begin{array}{cc}
A_{:T^\ell}^\top A_{:T^\ell} + \eta I&A_{\overline\Gamma_0^\ell T^\ell}^\top\\
A_{\overline\Gamma_0^\ell T^\ell}&I
\end{array}
\right]\left[
\begin{array}{cc}
\bfx_{T^\ell} \\
\bfy_{\overline\Gamma_0^\ell}
\end{array}
\right]=\left[
\begin{array}{r}
 A_{: T^\ell}^\top \epsilon{\bf 1} \\
 \epsilon {\bf 1}
\end{array}
\right].
\end{eqnarray*}
The solution,  $(\widetilde\bfx^\ell;\widetilde\bfy^\ell)$,  can be derived by
\begin{eqnarray*}
\widetilde\bfx^\ell_{T^\ell} &=& \left[
A_{\Gamma_0^\ell T^\ell}^\top A_{\Gamma_0^\ell T^\ell} + \eta I \right]^{-1}\left[
 A_{\Gamma_0^\ell T^\ell}^\top\epsilon {\bf 1}
\right],\\
\widetilde\bfx^\ell_{\overline T^\ell} &=&0,\\
\widetilde\bfy^\ell_{\Gamma_0^\ell} &=&0,\\
\widetilde\bfy^\ell_{\overline\Gamma_0^\ell} &=& \epsilon\bfe - A_{\overline\Gamma_0^\ell T^\ell}  \widetilde\bfx^\ell_{T^\ell}.
\end{eqnarray*}
If $\widetilde\bfy^\ell_{ \Gamma_-^\ell}\leq 0$, namely, $(\widetilde\bfx^\ell;\widetilde\bfy^\ell)$ is the solution to  (\ref{QP}), then we set $\bfv^\ell=(\widetilde\bfx^\ell;\widetilde\bfy^\ell)$. Otherwise, this point will not be taken into consideration,   and we set $\bfz^{\ell+1}=\bfu^\ell.$  The computational complexity of addressing the above equations is about $O(ms^2+s^3)$.

Overall, the computational complexity of each iteration is $$O(\sigma mn +ms^2+s^3).$$
 
\subsection{Convergence analysis}
The first result shows that Armijo-type step size {(see \cite{armijo1966minimization} for more details)} $\{\tau_\ell\}$ is well defined.
\begin{lemma}\label{alpha-well-defined} For any $0<\tau\leq \frac{1}{2\rho+2\lambda_{\max}}$, it holds that
		\begin{eqnarray}
		\label{zk-alpha-zk}f(\bfz^\ell(\tau))\leq f(\bfz^\ell) - \rho\|\bfz^\ell(\tau)-\bfz^\ell\|^2,
		\end{eqnarray}
and thus $\inf_{\ell\geq 0}\{\tau_\ell\}\geq \underline{\tau}>0$, where
				\begin{eqnarray}
		\label{lower-bd-alphak}
		\begin{array}{l}
		\underline{\tau}:=\min\left\{1,~\frac{\beta}{2\rho+2\lambda_{\max}}\right\}.
		\end{array}\end{eqnarray}
\end{lemma}

 \begin{lemma}\label{zk1-zk-0}  Let $\{\bfz^\ell\}$ be the sequence generated by \gpa\ and $\underline{\tau}$ be given by (\ref{lower-bd-alphak}). Then the following results hold.
 \begin{itemize}
 \item[i)] Sequence $\{\bfz^\ell\}$ is bounded and $\lim_{\ell\rightarrow\infty}\|\bfz^{\ell+1}-\bfz^\ell\|=\lim_{\ell\rightarrow\infty}\|\bfu^{\ell}-\bfz^\ell\|=0$.
 \item[ii)] Any accumulating point of $\{\bfz^\ell\}$ is a $\tau$-stationary point with $0<\tau\leq \underline{\tau}$ of problem  (\ref{ell-0-1}).
 \end{itemize}
\end{lemma}
 
The above lemma allows us to conclude that the whole sequence converges {to a unique local minimizer without any assumptions. But with an additional condition, the whole sequence can achieve a global minimizer.}
 \begin{theorem}\label{global-convergence}{The whole sequence, $\{\bfz^\ell\}$,  generated by \gpa\ converges to  a $\tau$-stationary point $\bfz^*$ of problem  (\ref{ell-0-1}), which is necessarily  
a unique local minimizer. If $\bfz^*$ is also a $\tau_*$-stationary point with $\tau_*\geq1/(2\lambda_{\min})$, where  $\tau_*$ is defined by \eqref{eta-point-global-tau*}, then it is a  global minimizer.}
\end{theorem}
{
The following theorem claims that \gpa\ can terminate at the limit of the sequence after a certain point.
 \begin{theorem}\label{convergence-rate}  Let $\{\bfz^\ell\}$ be the sequence generated by \gpa\ with $0<\rho\leq \lambda_{\min}$ and $\bfz^*$ be  its limiting point.   Then  \gpa\ will terminate at $\bfz^*$  within finite steps, namely, there is a finite $\kappa\geq 1$ such that
  \begin{eqnarray}
 \label{case-1-convergence}  \bfz^{\ell}=\bfv^{\kappa}=  \bfz^{*},~~\forall~\ell > \kappa.
\end{eqnarray} 
\end{theorem}}
{Finally, we can even show that the distance between every iterate $\bfx^\ell$ and true signal $\bfx^\tr$  has an upper bound.
 \begin{theorem}\label{convergence-rate-true} Let $\Phi$ and $\var$ be given as in Theorem \ref{oracle-property}, $\{\bfz^\ell\}$ be the sequence generated by \gpa\ and $\bfz^*$ be its limiting point. We choose $\eta$ as \eqref{eta-tuning} and $\rho$  as $0<\rho\leq \lambda_{\min}$. Then  there is a finite $\kappa\geq1$ such that, for any $\ell\geq \kappa$, \begin{eqnarray*}  
\mathbb{P} \{\| c_* {\bfx}^\ell-\bfx^\tr\| \leq 2{\rm sin}\left(\phi\pi\right) \} 
 \geq  1-e^{-2m\delta^2}-\theta.
  \end{eqnarray*} 
\end{theorem}}
\subsection{Comparisons with other methods}
{We would like to compare {\tt GPSP} with some other methods that have been proposed to deal with the noisy recovery with sign flips, namely, model \eqref{noise-1bcs}.    Their comparisons on the theoretical guarantees and the computational complexity are summarized in Table \ref{tab:compare-theory-complexity}. Only three methods {\tt OSL0}, {\tt PDASC} and {\tt GPSP} have been established convergence results.  Since {\tt PDASC} converges only when the initial point is chosen close to the accumulating point of the sequence, it has a local convergence property. Both  {\tt OSL0} and {\tt PDASC} also require some assumptions on data $\Phi$ or $A$ to derive the convergence. Moreover, we note that  {\tt WPDASC}, {\tt PDASC} and {\tt GPSP} are the second-order methods since they make use of Hessian matrix (i.e., the second-order information) of the objective function. Therefore, they have slightly higher computational complexity. Here, $s_\ell$ is the sparsity level of the point at $\ell$th iteration. However, empirical numerical experiments have demonstrated that the second-order methods can converge within much fewer steps and yield a much higher order of accuracy than the first-order methods. This is another reason that {\tt GPSP} is capable of delivering relatively desirable accuracy.
\begin{table}[!th]
	\renewcommand{\arraystretch}{1.25}\addtolength{\tabcolsep}{-1pt}
	\caption{Comparisons of different algorithms.}\vspace{-2mm}
	\label{tab:compare-theory-complexity}
	\begin{center}
		\begin{tabular}{lccr}
			\hline	
		Algs.  	&Convergence&Assumptions &Complexity\\
		  	& & on  $\Phi$ or $A$ &\\\hline
\multicolumn{4}{c}{First-order methods}\\\hline
{\tt RBIHT}	\cite{fu2014robust}& $--$	&  $--$ &	$O( mn )$\\
{\tt AOPF} \cite{yan2012robust}	& $--$	&  $--$  &$O( mn)$	\\
{\tt NARSS}	\cite{movahed2014recovering}& $--$	&  $--$ &	$O( mn )$\\
{\tt OSL0}	\cite{dai2016noisy}& Global	&  Yes  &	$O( mn )$\\
\hline
\multicolumn{4}{c}{Second-order methods}\\\hline
{\tt WPDASC}	\cite{fan2021robust}& $--$	&  $--$  &	$O( mn +ms_\ell^2+s_\ell^3)$\\
{\tt PDASC} \cite{huang2018robust}	& Local &Yes&$O( mn +ms_\ell^2+s_\ell^3)$\\
{\tt GPSP} & Global & No & $O(\sigma mn +ms^2+s^3)$\\	\hline
		\end{tabular}
	\end{center}
\end{table} }
{
\begin{remark}We now highlight the difference between  {\tt GPSP} and {\tt OSL0} in \cite{dai2016noisy} developed to solve \eqref{0-1-regular-approx}, which is a penalized version   our proposed model \eqref{ell-0-1}. Firstly, the  different frameworks  mean that {\tt GPSP} is a second-order method while {\tt OSL0} is a first-order method.  Moreover,   {\tt OSL0} converges to a local minimizer if the spectral norm of matrix $A$ is bounded by some chosen parameters and converges to a global minimizer if further assuming that the initial point is sufficiently close to this global minimizer. However, {\tt GPSP}  converges to a unique local minimizer without any assumptions on $A$ and converges to a global minimizer if $\tau_*\geq1/(2\lambda_{\min})$. Hence it does not impose any condition on the initial point. 
\end{remark}}

\section{Numerical Experiments}\label{sec:numerical}
\noindent In this section, we will conduct extensive numerical experiments  {to showcase the performance of our proposed} {\tt GPSP} (available at \url{https://github.com/ShenglongZhou/GPSP}), by using MATLAB (R2019a) on a laptop of  $32$GB memory and Inter(R) Core(TM) i9-9880H 2.3Ghz CPU. 

\subsection{Testing examples}
  Examples with the data generated from the Gaussian distributions are taken into account.
\begin{example}[Independent covariance \cite{yan2012robust, dai2016noisy}]\label{ex:cs-ind} Entries of $\Phi:=[\bfa_1,\cdots,\bfa_m]^\top\in\R^{m\times n}$ and the nonzero entries of ground-truth $s_*$-sparse vector $\bfx^\tr\in\R^n$ (i.e., $\|\bfx^\tr\|_0\leq s_*$)  are generated from the independent and identically distributed (i.i.d.) samples of the standard Gaussian distribution, $\mathcal{N}(0,1)$. To avoid tiny nonzero entries of $\bfx^\tr$, let $x^\tr_i= x^\tr_i+{\rm sign}(x^\tr_i)$ for nonzero $x^\tr_i$, followed by normalizing  $\bfx^\tr$ to be a unit vector. Let $\bfc^\tr={\rm sgn}(\Phi\bfx^\tr)$ and $\bfc ={\rm Diag}(\bfh) {\rm sgn}(\Phi\bfx^\tr+ \beps)$, where entries of noise $\beps$  are the i.i.d. samples of $\mathcal{N}(0,0.1^2)$ and  $\lceil r m\rceil$ entries of $\bfh $ are randomly selected to be $-1$, where  $r$  is the flipping ratio.
\end{example}

\begin{example}[Correlated covariance \cite{huang2018robust}]\label{ex:cs-cor} Rows of $\Phi$ are generated from the i.i.d. samples  of $\mathcal{N}(0,\Sigma)$ with $\Sigma_{ij}=v^{|i-j|}, i,j\in[n]$, where $v\in(0,1)$. Then $\bfx^\tr, \bfc^\tr$, and $\bfc$ are generated the same as those in \Cref{ex:cs-ind}.
\end{example}
To demonstrate the performance of one method, apart from the  CPU {\tt TIME}, we will also report the signal-to-noise ratio ({\tt SNR}) in dB, the Hamming error ({\tt HE}) and the Hamming distance ({\tt HD}). They are defined by
\begin{eqnarray*}
{\tt SNR} &:=&-20{\log}_{10} \| \bfx -\bfx^\tr\|,\\
{\tt HD}&:=&({1}/{m}) \|{\rm sgn}(\Phi\bfx)-\bfc\|_0,\\
{\tt HE}&:=& ({1}/{m})  \|{\rm sgn}(\Phi\bfx)-\bfc^\tr\|_0,
\end{eqnarray*}
where $\bfx$ is the solution obtained by one method. The larger {\tt SNR} (or the smaller  {\tt HE} or {\tt HD}) means the better recovery.

\subsection{Implementation and parameter  selection}
We terminate {\tt GPSP} if $\|\bfz^{\ell}-\bfu^\ell\|\leq \varepsilon$ or $\ell>2000$ and set $\beta=0.5$, $\rho=10^{-6}, \varepsilon=10^{-4}$ and $\bfz^0=0$ in all experiments. Parameters $\epsilon$, $\eta$, $s$ and $k$ in \eqref{ell-0-1} are tuned as follows. 

\begin{figure}[!th]
\centering
\begin{subfigure}{.24\textwidth}
	\centering
	\includegraphics[width=1.03\linewidth]{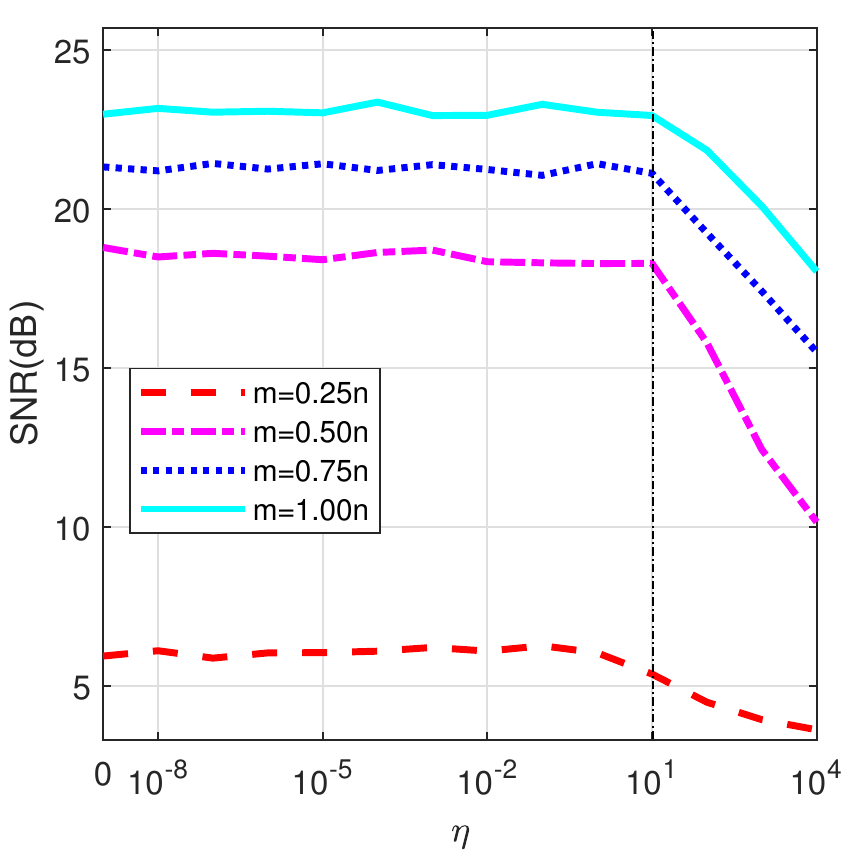}
\end{subfigure}
\begin{subfigure}{.24\textwidth}
	\centering
	\includegraphics[width=1.03\linewidth]{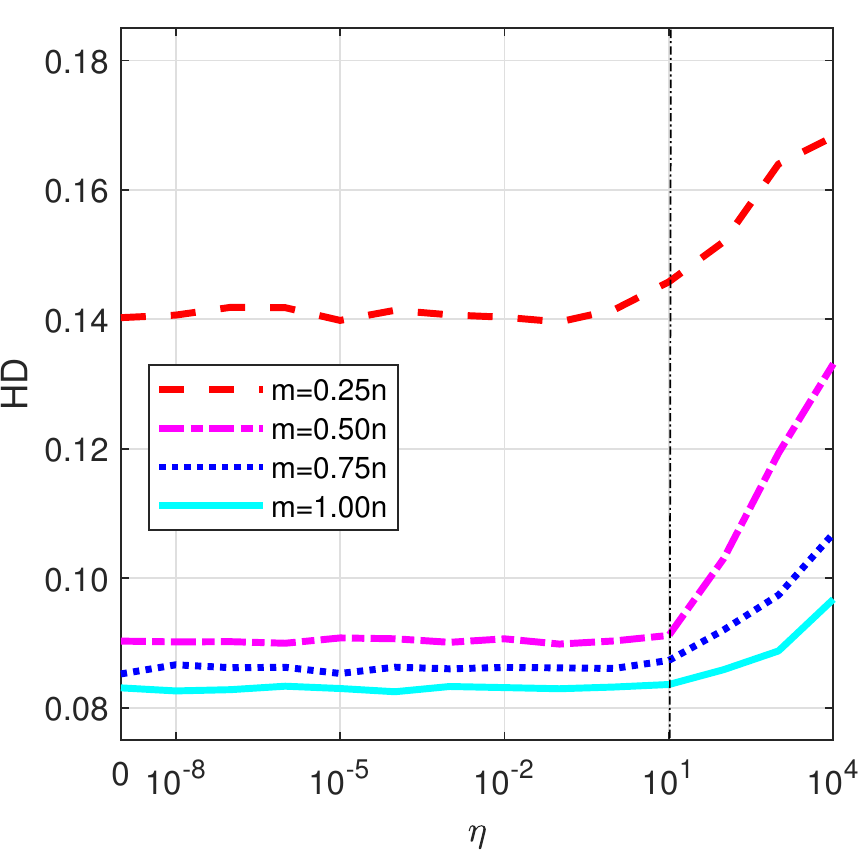}
\end{subfigure} 
\caption{Effect of $\eta$ for  \Cref{ex:cs-ind}.}
\label{fig:ex1-cs-eta}  
\end{figure}

{\bf (i) Selection of $\eta$.} {Although Theorem \ref{oracle-property} suggested that $\eta$ should be set as \eqref{eta-tuning}, it is difficult to estimate the value of the right hand side of \eqref{eta-tuning}. Therefore, we tested a wide range of $\eta$ (e.g, $\eta\in[0,10^4]$) to see how it effects our algorithm.} To proceed with that, for \Cref{ex:cs-ind}, we fix $(n,s_*,r,\epsilon)=(500,5,0.05,0.01)$ and  $k=\lceil 0.01m\rceil$  but vary $m\in\{0.25,0.5,0.75,1\}n$ and $\eta\in\{0, 10^{-8},10^{-7},\cdots,10^4\}$. 
Average results over 200  {trials} are reported in \Cref{fig:ex1-cs-eta}. It can be evidently seen that results are stabilized when $\eta\in[0, 10]$ while getting worse when $\eta>10$ is rising. Similar trends are also  observed for \gpa\ solving \Cref{ex:cs-cor}.  This well testifies that $\eta$ should be chosen smaller than a threshold, as shown in \eqref{eta-tuning}.  Therefore, any value in $[0,10]$ can be used to set $\eta$.   For simplicity, we fix $\eta=10^{-4}$.

\begin{figure}[!th]
\centering
\begin{subfigure}{.24\textwidth}
	\centering
	\includegraphics[width=1.03\linewidth]{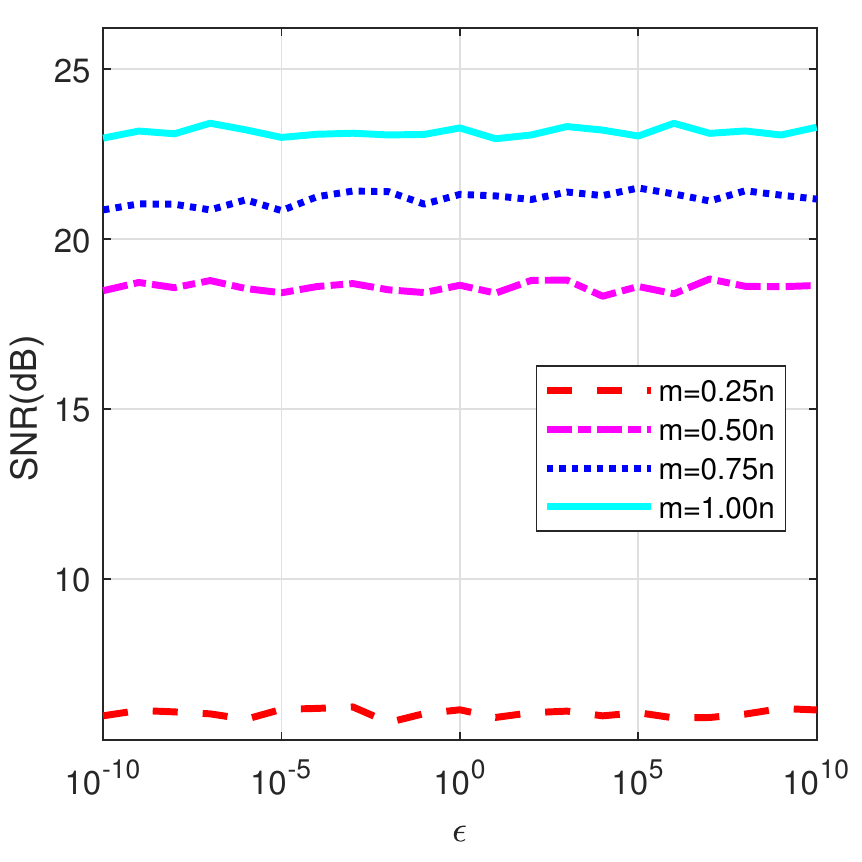}
\end{subfigure} 
\begin{subfigure}{.24\textwidth}
	\centering
	\includegraphics[width=1.03\linewidth]{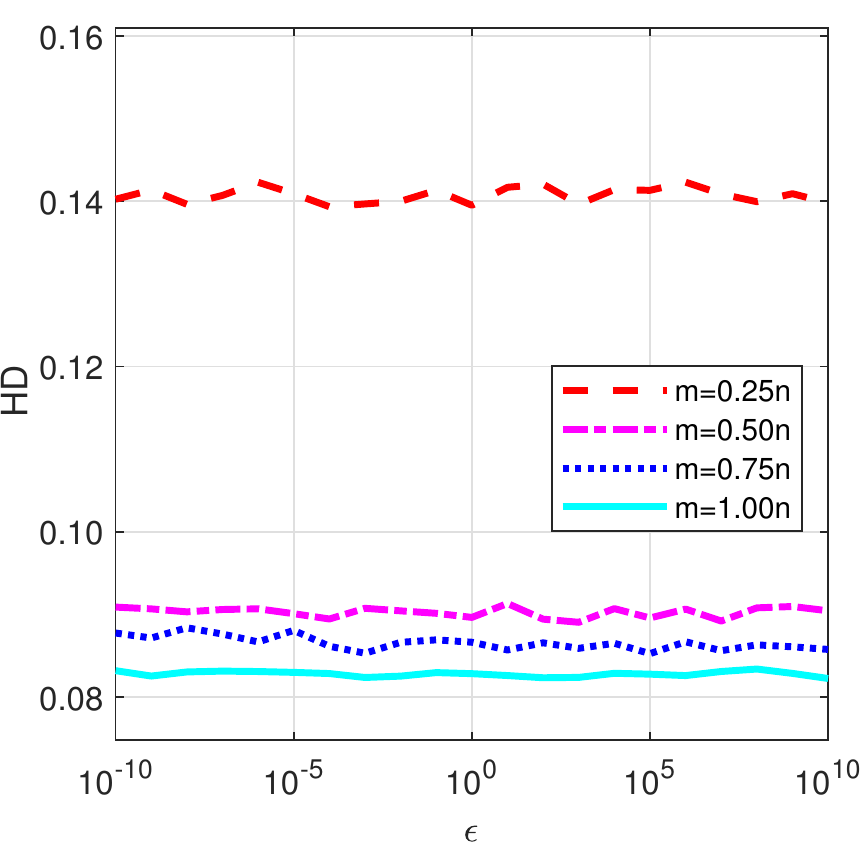}
\end{subfigure} 
\caption{Effect of $\epsilon$ for  \Cref{ex:cs-ind}.}
\label{fig:ex1-cs-eps}  
\end{figure}

{\bf (ii) Selection of $\epsilon$.} {As mentioned in Remark \ref{remark-para}, {\tt GPSP} might be quite robust to the choices of $\epsilon$ since the error bound theory (i.e., Theorem \ref{oracle-property}) does not impose any assumptions on  $\epsilon$. To testify that,  we fix  $(n,s_*,r)=(500,5,0.05)$ and  $k=\lceil 0.01m\rceil$  but alter  $m\in\{0.25,0.5,0.75,1\}n$ and $\epsilon\in[10^{-10},10^{10}]$ for \gpa\   solving \Cref{ex:cs-ind}. As presented in \Cref{fig:ex1-cs-eps}, the  average results over 200  {trials} do not fluctuate significantly when varying $\epsilon$, which shows  robustness of \gpa\ to the choices of $\epsilon$ in range  $[10^{-10},10^{10}]$. An underlying explanation is as follows: In our numerical experiments, we observed that $\|\bfx^\ell\|$ was increasing along with the rising of $\epsilon$. So the final iteration, $\bfx^\ell$,  is indeed influenced by  $\epsilon$. However, after the normalization, namely, $\overline{\bfx}=\bfx^\ell/\|\bfx^\ell\|$, the impact of $\epsilon$ is eliminated.} Since similar performance can be seen for \gpa\ solving \Cref{ex:cs-cor}, the corresponding results are omitted.   For simplicity, we fix $\epsilon=0.01$ in the subsequent numerical experiments.

\begin{figure}[!th]
\centering
\begin{subfigure}{.24\textwidth}
	\centering
	\includegraphics[width=1\linewidth]{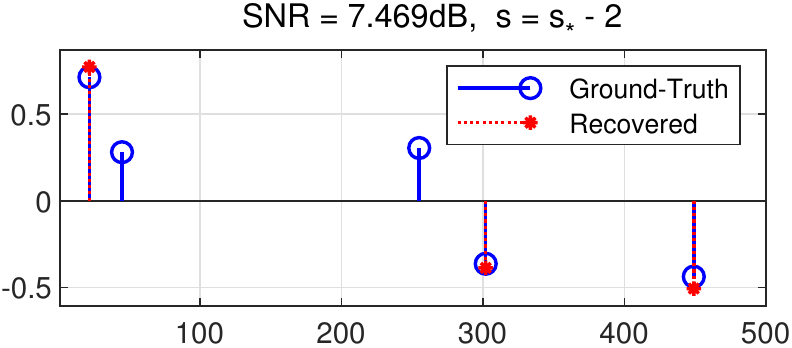}
\end{subfigure}
\begin{subfigure}{.24\textwidth}
	\centering
	\includegraphics[width=1\linewidth]{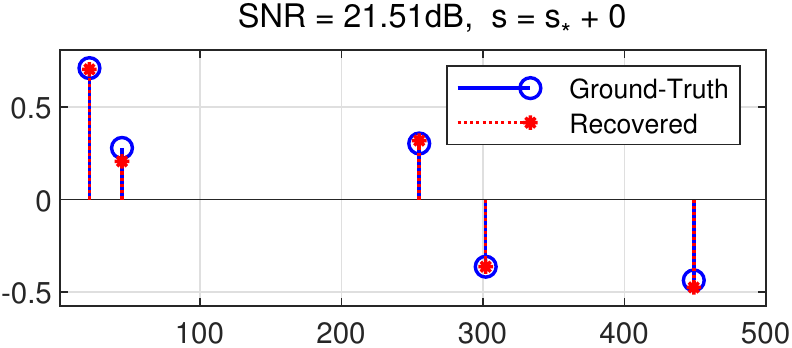}
\end{subfigure} \\\vspace{2mm}
\begin{subfigure}{.24\textwidth}
	\centering
	\includegraphics[width=1\linewidth]{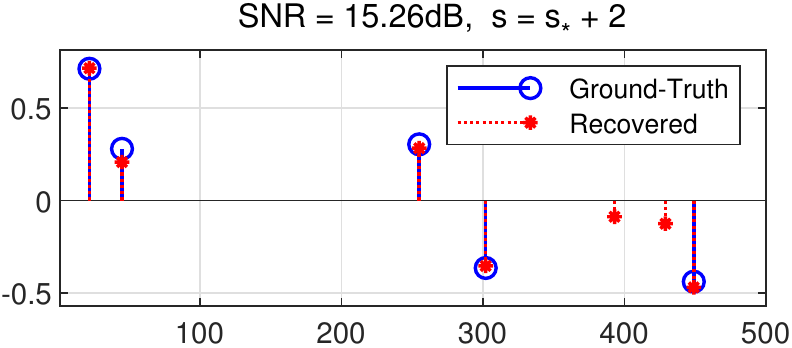}
\end{subfigure}
\begin{subfigure}{.24\textwidth}
	\centering
	\includegraphics[width=1\linewidth]{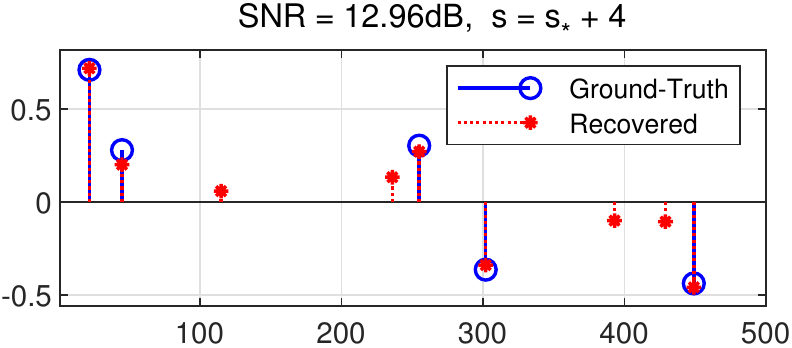}
\end{subfigure}
\caption{Effect of $s$ for  \Cref{ex:cs-ind}.}
\label{fig:ex1-cs-s}  
\end{figure}

{\bf (iii) Selection of $s$.} Sparsity level $s$ clearly has a heavy influence of  {the recovery quality.} 
As shown in \Cref{fig:ex1-cs-s} where $(n,m,s_*,r)=( 500,250,5,0.05)$, the ground-truth signal, $\bfx^\tr$, has $s_*=5$ nonzero components with their indices denoted by $T^\tr=\supp(\bfx^\tr)$. Apparently, \gpa\ gets the most accurate signal if we set $s=s_*$ because it almost exactly recovers those nonzero components. For  $s=s_*-2=3$, the recovered signal has  $3$ nonzero components, however, their indices belong to the true support set, $T^\tr$.
While for $s=s_*+2$ or $s=s_*+4$,   \gpa\ generates a solution  $\bfx$ whose support set $\supp(\bfx)$ covers $T^\tr$  {with extra} 
incorrect indices. However, compared with magnitude $|x_i|, i\in  T^\tr$, those  redundant nonzero components $|x_i|, i\in \supp(\bfx)\setminus T^\tr$ are pretty small. If we remove those small parts and normalize the signal to have a unit length, then the new signal is much closer to $\bfx^\tr$. For simplicity, we set $s=s_*$ in the sequel.

\begin{figure}[!th]
\centering
\includegraphics[width=1\linewidth]{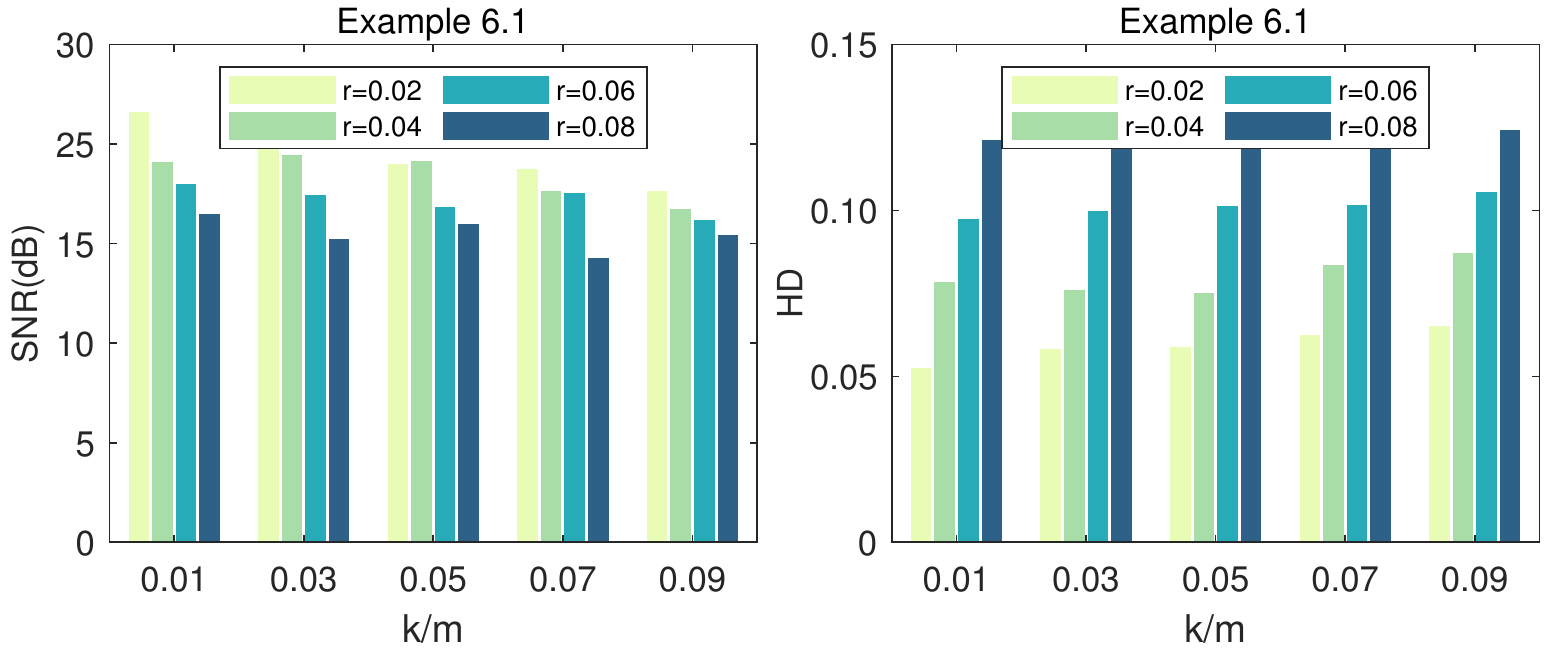}\\\vspace{2mm}
\includegraphics[width=1\linewidth]{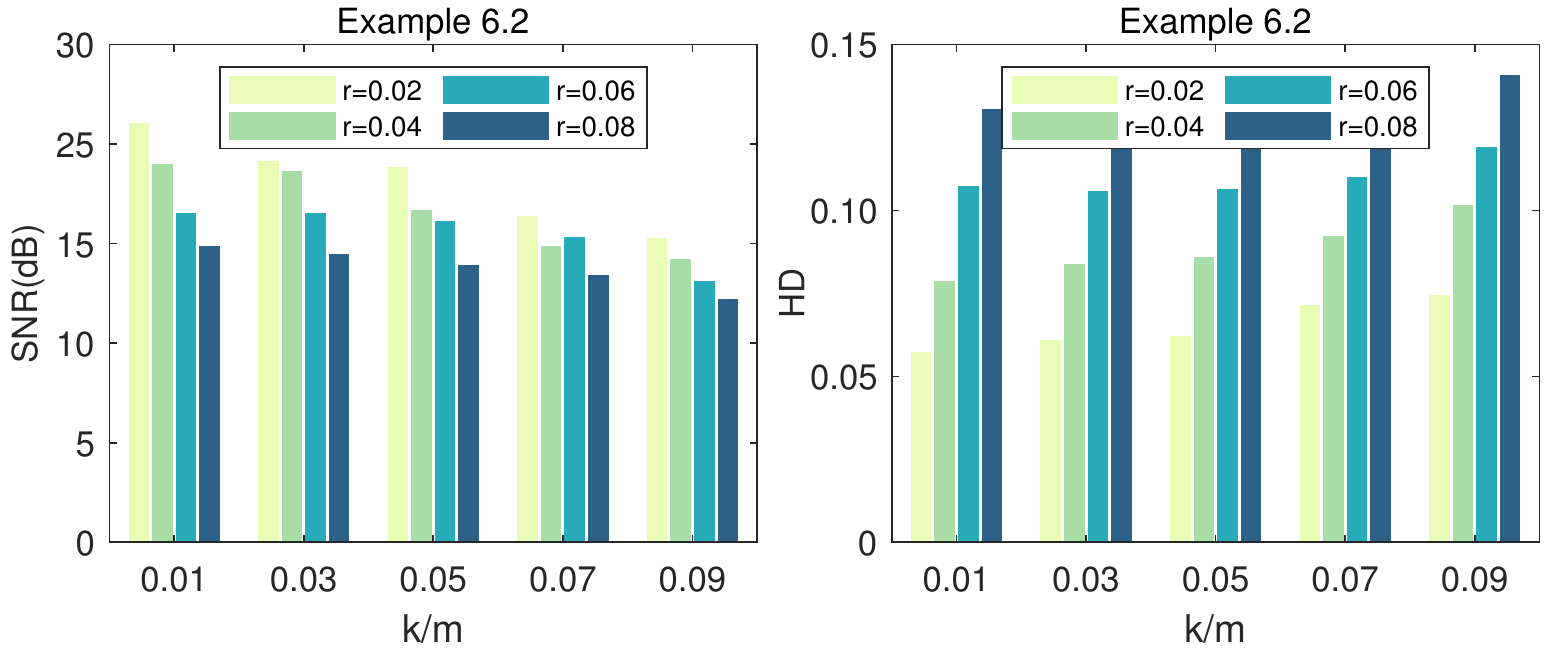}
\caption{Effect of $k$ for  \Cref{ex:cs-ind} and  \Cref{ex:cs-cor}.}
\label{fig:ex1-cs-k} \vspace{-3mm}
\end{figure}

{\bf (iv) Selection of $k$.}  Note that $k$ is the upper bound of the number of sign flips  of $\Phi\bfx$ and is usually unknown beforehand. However, model \eqref{ell-0-1} does not require an exact $k$. One could either fix it by a small integer (e.g. $k=\lceil 0.01m\rceil$) or start with a slightly bigger value  and reduce it iteratively. We tested \gpa\ for solving \Cref{ex:cs-ind} and \Cref{ex:cs-cor} under both schemes and corresponding numerical performance does not have a big difference. For instance, as indicated in \Cref{fig:ex1-cs-k}, where $(n,m,s_*,v)=(500,250,5,0.5)$, we select $k/m\in\{0.01,0.03,0.05,0.07,0.09\}$ and then fix it for {\tt GPSP}. Evidently, for each case of flipping ratio $r$,  results {\tt SNR} and {\tt HD} do not vary significantly along with $k$ altering.

 {We note that in the above numerical experiments, the upper bound, $k$, was  unequal to true number of sign flips $rm$ but close to $rm$.   
However, we have also tested much larger $k/m$ (e.g., 0.2) than $r$ (e.g., 0.02), the recovery accuracies were degraded greatly. Therefore, we could conclude that \gpa\ is quite robust to these $k$  around the true number of sign flips,   which suggests that $k$ should not be chosen too far away from $rm$.}
Hence, in our numerical experiments, we pick $k=\lceil 0.01m\rceil$ if no additional information is provided.

\subsection{Benchmark methods} {Six state-of-the-art solvers are selected for comparisons. They are  {\tt BIHT} \cite{jacques2013robust}, {\tt AOPF} 
 ({\tt BIHT-AOP-flip}, \cite{yan2012robust}), {\tt PAOPF} 
 ({\tt PIHT-AOP-flip}, \cite{huang2018pinball}), {\tt PIHT} (This code was implemented based on the pinball loss from \cite{huang2018pinball} to address the noisy recovery with sign flips),  {\tt PDASC}~\cite{huang2018robust} and {\tt WPDASC} \cite{fan2021robust}.  We note that  {\tt BIHT} behaves well for noiseless recovery but has shown undesirable performance for noisy recovery. We include it in comparison  since we use its obtained solution $\bfx$ to generate the upper bound of the number of sign flips for {\tt AOPF}, {\tt PAOPF}, and {\tt PIHT}, namely, set $L=\|{\rm sgn}(\Phi\bfx)-\bfc\|_0$. This is one option to generate the upper bound  suggested by \cite{yan2012robust}.} Like our method, the first four methods need to specify sparsity level $s$ that is set by $s=s_*$.    Moreover, to accelerate the termination of {\tt PIHT}, we stop it when the number of wrong recovery signs is smaller than $L$ and the number of iterations is over 50. The other parameters for each method are chosen to be their default values. All methods are initialized by $\bfx^0= 0$, and their final solutions are normalized to have a unit length. Finally, to make fair comparisons,  we keep the first $s$ largest absolute values for {\tt PDASC} and  {\tt WPDASC} if one obtains a solution with more than $s$ nonzero entries. Overall, all methods start with the same initial points and generate solutions with at most $s$ nonzero entries.

\subsection{Numerical comparisons}
 We now apply the seven methods into solving  two examples under different scenarios. For each scenario, we report average results over $200$ instances if $n\leq1000$ and $20$ instances otherwise. Note that for each example, there are five factors $(n,m,s_*,r,v)$, where $v$ only makes sense for \Cref{ex:cs-cor}. In the following numerical comparisons, we shall see the effect of these factors by altering one factor while fixing the others.

\begin{figure*}[!th]
\centering 
\begin{subfigure}{1\textwidth}
	\centering
	\includegraphics[width=.99\linewidth]{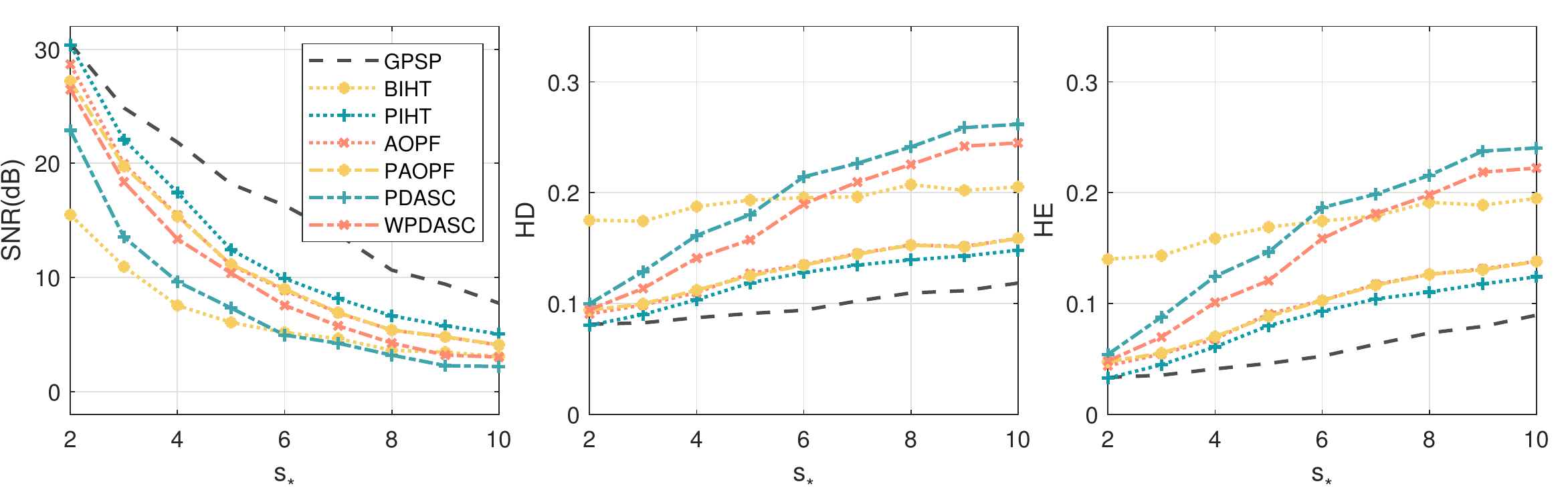}
	\caption{Effect of $s_*$ for \Cref{ex:cs-ind}.}
	\label{fig:ex1-2-cs-s}
\end{subfigure}	\\\vspace{-1mm}
 
 \begin{subfigure}{1\textwidth}
	\centering
	\includegraphics[width=.99\linewidth]{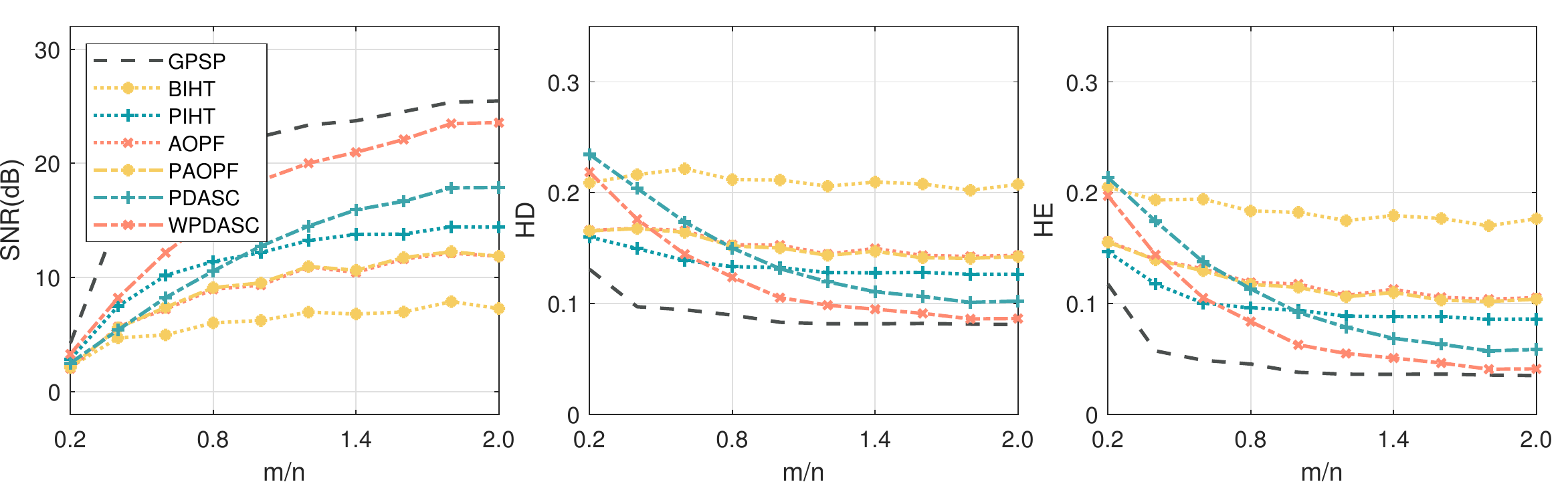}
	\caption{Effect of $m$ for \Cref{ex:cs-cor}.}
	\label{fig:ex1-2-cs-m}
\end{subfigure}	\\\vspace{-1mm}

 \begin{subfigure}{1\textwidth}
	\centering
	\includegraphics[width=.99\linewidth]{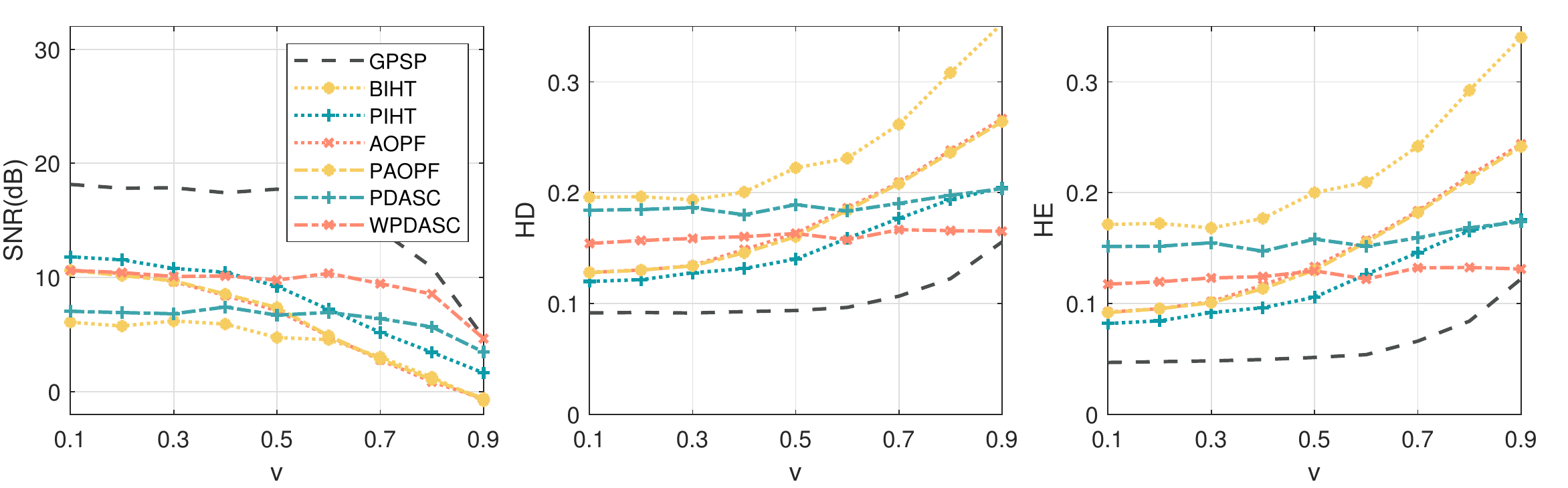}
	\caption{Effect of $v$ for \Cref{ex:cs-cor}.}
	\label{fig:ex1-2-cs-v}
\end{subfigure}	 
\caption{Effect of $(s_*,m, v)$ for seven algorithms.}
\label{fig:ex1-2-cs-all} 
\end{figure*}

   {\bf (a) Effect of $s_*$.}
We first employ seven methods to solve \Cref{ex:cs-ind} and increase $s_*$ from $2$ to $10$ with fixing  $(n,m,r,v)=(500,250,0.05,0.5)$. As shown in \Cref{fig:ex1-2-cs-s},   {\tt GPSP} gets the highest {\tt SNR}, the smallest {\tt HD} and {\tt HE} for each $s_*$, followed by {\tt AOP}, {\tt PAOP} and {\tt PIHT}.   
The lines of  {\tt SNR}  display declining trends, which means the signal is getting harder to recover when it has more  nonzero components, namely, $s_*$ is getting bigger.

 {\bf (b) Effect of $m$.} To see the effect of sample size $m$, 
we  select it from range $\{0.1,0.3,\cdots,1.5\}n$ and fix $(n,s_*,r,v)=(500, 5, 0.05,0.5)$. As shown in \Cref{fig:ex1-2-cs-m}, again, {\tt GPSP} outperforms the others for solving \Cref{ex:cs-cor} since it delivers much higher {\tt SNR} and lower {\tt HD} and {\tt HE}.  It is evidently seen that all methods are behaving better along with the rising of sample size $m$ because the signal is getting easier to recover when more samples are available. 

{\bf (c) Effect of $v$.} We note that in \Cref{ex:cs-cor}, the larger $v$ is, the more correlated each pair of samples (i.e., rows in $\Phi$) are, leading to more difficult recovery.  To see this,  we alter $v$ from $\{0.1,0.2,\cdots,0.9\}$ but fix  $(n,m,s_*,r)=(500, 250,5,0.05)$, and report the average results in \Cref{fig:ex1-2-cs-v}. As expected, the larger $v$ is, the more difficult the recovery is.  It is observed that {\tt GPSP} is quite robust to $v$ between $0.1$ and $0.6$ since the produced results  stay steadily when  $v\in[0.1,0.6]$. No matter how $v$ changes,  {\tt GPSP} always performs the best results among those  methods.

\begin{figure*}[!th]
\centering
\includegraphics[width=1\linewidth]{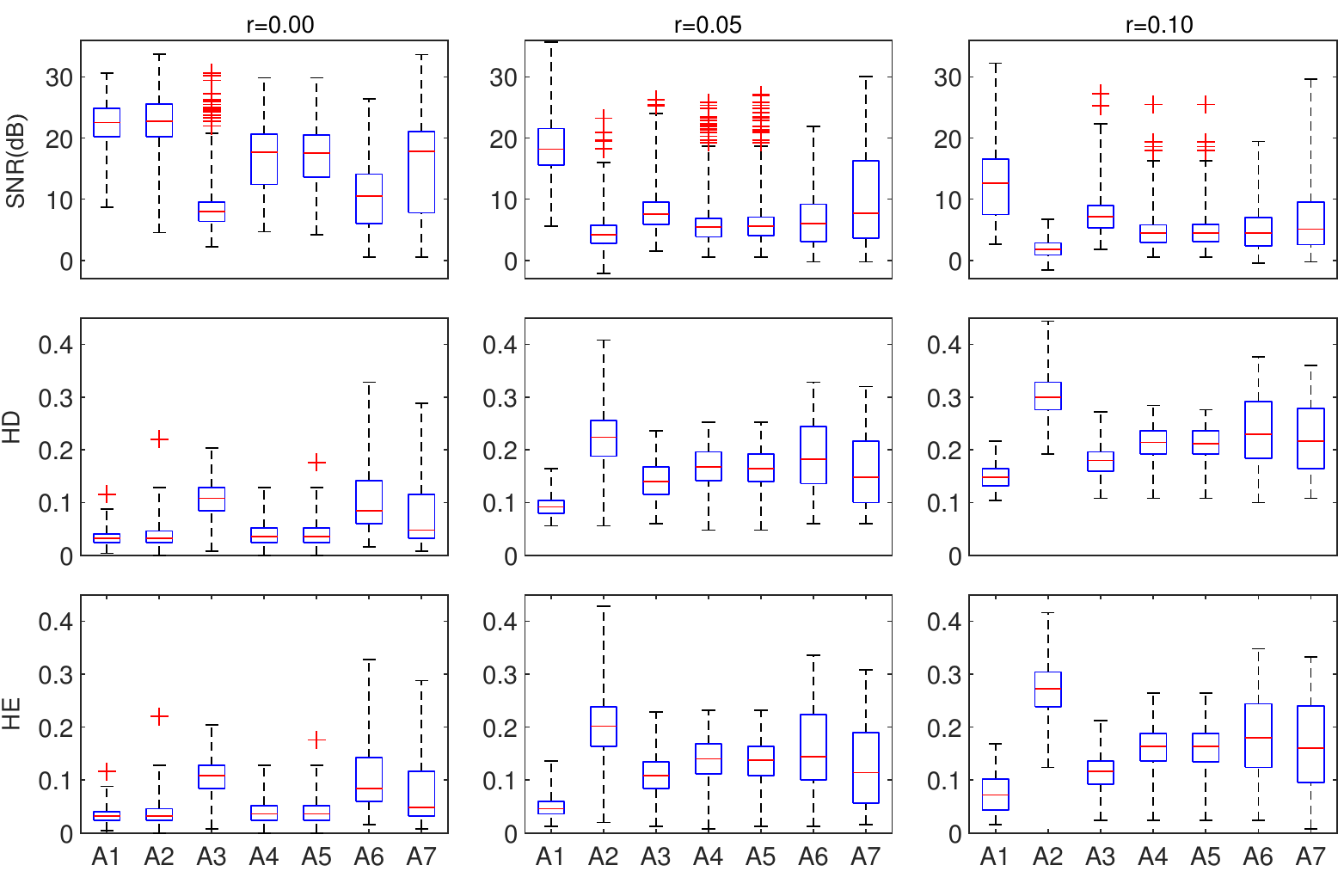}
\caption{Effect of $r$ for \Cref{ex:cs-cor}, where A1-A7 stand for {\tt GPSP}, {\tt BIHT}, {\tt PIHT}, {\tt AOPF},  {\tt PAOPF}, {\tt PDASC},   {\tt WPDASC}, respectively.}
\label{fig:ex1-2-cs-r}   
\end{figure*}

{\bf (d) Effect of $r$.} To see the effect of flipping ratio $r$, we alter it from $\{0,0.05,0.1\}$ but fix  $(n,m,s_*,v)=(500, 250,5,0.5)$.   The box-plots of each method for solving \Cref{ex:cs-cor}  are presented in \Cref{fig:ex1-2-cs-r}. In each box, the central mark (red line) indicates the median,  the bottom, and top edges of the box indicate the $25$th and $75$th percentiles, respectively. The outliers are plotted individually using the `+' symbols.

 As expected, the larger $r$ is, the worse performance of each method, because more correct signs are flipped. This can be testified by {\tt SNR} (resp. {\tt HD} and {\tt HE}) whose median obtained by each method is declining (resp. rising) when $r$ ascends. Once again,   {\tt GPSP} behaves the best because it delivers the highest median of {\tt SNR} and the lowest median of {\tt HD} and {\tt HE} in each box. Similar results can be observed for \Cref{ex:cs-ind}  and are omitted here.
 
\begin{table}[!th]
	\renewcommand{\arraystretch}{1}\addtolength{\tabcolsep}{-3pt}
	\caption{Effect of the bigger values of $n$.}\vspace{-2mm}
	\label{fig:ex2-cs-ind-n}
	\begin{center}
		\begin{tabular}{ll ccccccc}
			\hline			
	&	$n$	&	{\tt GPSP} & {\tt BIHT}	&	{\tt PIHT} &	{\tt AOPF}	&	{\tt PAOPF}		& 	{\tt PDASC}	&	{\tt WPDASC}		\\\hline
	&&\multicolumn{7}{c}{ \Cref{ex:cs-ind}}\\\cline{3-9}
{\tt SNR}	&	$5000$	&	15.51 	&	4.612 	&	9.584 	&	7.948 	&	7.948 	&	-0.228 	&	0.786 	\\
(dB)	&	$10000$	&	12.17 	&	4.431 	&	9.230 	&	7.464 	&	7.464 	&	-1.755 	&	-1.141 	\\
	&	$15000$	&	12.47 	&	4.720 	&	8.976 	&	7.256 	&	7.256 	&	-2.023 	&	-1.477 	\\
	&	$20000$	&	12.89 	&	4.693 	&	9.095 	&	7.348 	&	7.348 	&	-2.135 	&	-2.111 	\\\hline
{\tt HD}	&	$5000$	&	0.092 	&	0.201 	&	0.122 	&	0.138 	&	0.138 	&	0.353 	&	0.332 	\\
	&	$10000$	&	0.106 	&	0.207 	&	0.125 	&	0.141 	&	0.141 	&	0.425 	&	0.394 	\\
	&	$15000$	&	0.103 	&	0.202 	&	0.127 	&	0.142 	&	0.142 	&	0.438 	&	0.412 	\\
	&	$20000$	&	0.102 	&	0.201 	&	0.126 	&	0.143 	&	0.143 	&	0.444 	&	0.443 	\\\hline
{\tt HE}	&	$5000$	&	0.051 	&	0.180 	&	0.087 	&	0.106 	&	0.106 	&	0.339 	&	0.317 	\\
	&	$10000$	&	0.067 	&	0.185 	&	0.089 	&	0.109 	&	0.109 	&	0.417 	&	0.383 	\\
	&	$15000$	&	0.065 	&	0.180 	&	0.092 	&	0.111 	&	0.111 	&	0.433 	&	0.404 	\\
	&	$20000$	&	0.062 	&	0.178 	&	0.091 	&	0.111 	&	0.111 	&	0.437 	&	0.437 	\\\hline
{\tt TIME}	&	$5000$	&	0.147 	&	0.528 	&	1.187 	&	0.072 	&	0.014 	&	0.378 	&	0.473 	\\
(s)	&	$10000$	&	0.580 	&	4.575 	&	6.124 	&	0.258 	&	0.056 	&	1.451 	&	1.698 	\\
	&	$15000$	&	1.054 	&	10.09 	&	14.39 	&	0.558 	&	0.139 	&	2.973 	&	3.791 	\\
	&	$20000$	&	1.807 	&	17.78 	&	24.82 	&	0.978 	&	0.243 	&	4.788 	&	6.459 	\\\hline

		&&\multicolumn{7}{c}{ \Cref{ex:cs-cor}}\\\cline{3-9}
{\tt SNR}	&	$5000$	&	13.35 	&	6.865 	&	7.368 	&	4.794 	&	5.351 	&	-0.948 	&	-0.860 	\\
(dB)	&	$10000$	&	11.55 	&	7.699 	&	7.180 	&	4.980 	&	5.630 	&	-1.857 	&	-1.340 	\\
	&	$15000$	&	11.22 	&	8.183 	&	6.982 	&	4.668 	&	5.386 	&	-1.802 	&	-1.896 	\\
	&	$20000$	&	11.67 	&	8.185 	&	6.942 	&	4.774 	&	5.549 	&	-2.224 	&	-2.324 	\\\hline
{\tt HD}	&	$5000$	&	0.099 	&	0.200 	&	0.148 	&	0.171 	&	0.168 	&	0.379 	&	0.375 	\\
	&	$10000$	&	0.106 	&	0.193 	&	0.148 	&	0.164 	&	0.163 	&	0.429 	&	0.405 	\\
	&	$15000$	&	0.109 	&	0.195 	&	0.149 	&	0.165 	&	0.163 	&	0.426 	&	0.431 	\\
	&	$20000$	&	0.106 	&	0.191 	&	0.150 	&	0.162 	&	0.162 	&	0.446 	&	0.452 	\\\hline
{\tt HE}	&	$5000$	&	0.058 	&	0.170 	&	0.113 	&	0.143 	&	0.138 	&	0.366 	&	0.362 	\\
	&	$10000$	&	0.070 	&	0.164 	&	0.115 	&	0.137 	&	0.133 	&	0.422 	&	0.396 	\\
	&	$15000$	&	0.072 	&	0.164 	&	0.117 	&	0.139 	&	0.135 	&	0.417 	&	0.423 	\\
	&	$20000$	&	0.069 	&	0.162 	&	0.117 	&	0.136 	&	0.133 	&	0.441 	&	0.448 	\\\hline
{\tt TIME}	&	$5000$	&	0.155 	&	0.542 	&	1.274 	&	2.201 	&	0.086 	&	0.386 	&	0.474 	\\
(s)	&	$10000$	&	0.566 	&	4.523 	&	6.448 	&	13.07 	&	1.600 	&	1.468 	&	1.736 	\\
	&	$15000$	&	1.139 	&	10.55 	&	14.71 	&	31.91 	&	3.575 	&	2.947 	&	3.726 	\\
	&	$20000$	&	1.996 	&	18.39 	&	25.44 	&	58.05 	&	10.60 	&	5.103 	&	7.150 	\\\hline

		\end{tabular}
	\end{center}
\end{table}

{\bf (e) Effect of $n$.} To see the computational speed of each method, we consider some bigger values of $n$ from $\{5000,10000,15000,20000\}$ with fixing  $(m,s_*,r,v)=(n/2, n/100,0.05,0.5)$. The average results are recorded in  \Cref{fig:ex2-cs-ind-n}. Obviously, {\tt GPSP}  {achieves the highest recovery accuracy in terms of the highest {\tt SNR}, the lowest {\tt HD} and {\tt HE}, against the other methods}. For the computational speed, {\tt PAOPF} runs the fastest for \Cref{ex:cs-ind} while {\tt GPSP} is the winner for   \Cref{ex:cs-cor}.

{In summary, we can conclude that {\tt GPSP} always achieved the highest recovery accuracies and ran the fastest among all benchmark methods. There are several reasons for its high performance. Firstly, our optimization model \eqref{ell-0-1} makes use of both sparsity and the number of sign flips and hence improves the recovery accuracy greatly. Moreover, as a second-order method, it is capable of enhancing the accuracy once again. Finally, it has a very low computational complexity and has been proven to terminate with finitely many steps, see \eqref{case-1-convergence}, thereby running super fast.
}

 \section{Conclusion}\label{sec:conclude}
\noindent  In this paper, we have proposed  {a new approach, {\tt GPSP}, that was derived from a non-convex formulation, \eqref{ell-0-1}, of the one-bit CS problem. The formulation has a simple objective function and exploits the double-sparsity to constrain the sparsity of the signal and the number of sign flips. Moreover, {\tt GPSP} can be deemed as a second-order method. Thanks to these two aspects, {\tt GPSP} is capable of delivering a high reconstruction quality.} To conquer the hardness resulting from the nonconvex and discrete constraints, we have established necessary and sufficient optimality conditions via the so-called $\tau$-stationarity. These optimality conditions have facilitated the design of a gradient projection subspace pursuit method {\tt GPSP}, which has been shown to admit global convergence and highly efficient numerical performance both in computation time and the recovery accuracy of signals.

 {There are a few topics pointed out by our referees that are worth exploring further. For instance,   fairly recently, authors in \cite{gianelli2016one, zhu2019computationally, ren2019sinusoidal} succeeded in extending the noisy one-bit CS to the scenario with a time-varying threshold being given as a prior condition for the sake of further improving the recovery accuracy.   Therefore, it is interesting to apply the developed techniques and the proposed algorithm into such a case. 
 }


%

\appendices
\section{Proofs of all theorems in Section \ref{sec:opt}} \label{app:opt}

\subsection{Proof of \Cref{nec-suff-opt-con-KKT} }
\begin{IEEEproof}
{\bf a) Necessity.} Let $\bfz^*$  be  a local  minimizer  of (\ref{ell-0-1}). Then there is a $\delta>0$ such that, for any $\bfz\in \F \cap N(\bfz^*, \delta)$,
\begin{eqnarray} \label{local-opt} 0&\leq& f(\bfz) - f(\bfz^*)\\
& \overset{\eqref{ell-0-1-obj-convex-smooth}}{\leq}&
 \langle \nabla f(\bfz^*), \bfz-\bfz^*\rangle +  \lambda_{\max} \|  \bfz-\bfz^*\|^2 =: g(\bfz).\nonumber
  \end{eqnarray}
 {To verify that $\bfz^*$ satisfies \eqref{KKT-point}, we consider the four cases.}

$\bullet$ $\|\bfy^*_+\|_0<k$. If there is an $i\in[m]$ such that  $(\nabla_{\bfy} f(\bfz^*))_i\neq 0$, then  {for any $t\in\R$, define}
\begin{eqnarray} 
\bfx_t&:=&\bfx^*+t\cdot 0,\nonumber\\
\label{x'-y'}  \bfy_t&:=&\bfy^*+t(\nabla_{\bfy} f(\bfz^*))_i \cdot{\bf e}_i,\\
\bfz_t&:=&(\bfx_t;\bfy_t).\nonumber 
\end{eqnarray}
 where  ${\bf e}_i \in\R^m$ is the $i$th column of the identity matrix. It is easy to see that $\bfz_t\in\F$ and
 $$g(\bfz_t)=(t+\lambda_{\max}t^2)[(\nabla_{\bfy} f(\bfz^*))_i]^2.$$
 {For any $t\in \left(\max\left\{-\delta/\left|(\nabla_{\bfy} f(\bfz^*))_i\right|, -1/\lambda_{\max}\right\},0\right)$, one can verify that}
$\bfz_t\in\F \cap N(\bfz^*, \delta)$ and $g(\bfz_t)<0$,
contradicting with \eqref{local-opt}. Thus $\nabla_{\bfy} f(\bfz^*)= 0$.

$\bullet$ $\|\bfy^*_+\|_0=k$.  If there is an $i\in\Gamma^*$ such that $(\nabla_{\bfy} f(\bfz^*))_i\neq 0$, then let $\bfx_t, \bfy_t$ and $\bfz_t$ be as \eqref{x'-y'}. The same reasoning for case  $\|\bfy^*_+\|_0<s$ enables to prove $(\nabla_{\bfy} f(\bfz^*))_i= 0$. This displays $(\nabla_{\bfy} f(\bfz^*))_{\Gamma^*} = 0$.  {To show   $(\nabla_{\bfy} f(\bfz^*))_{\overline\Gamma^*} \leq 0$, consider any given $i\in \overline\Gamma^*$. For any $t\leq 0$, let}
 \begin{eqnarray} \label{x'-y'-1}
\bfx_t:=\bfx^*+t\cdot 0,~~  \bfy_t:=\bfy^*+t \cdot{\bf e}_i, ~~\bfz_t:=(\bfx_t;\bfy_t).
\end{eqnarray}
It follows from $y_i^*=0$ and $t\leq 0$ that $\left(\bfy_t\right)_+=\bfy^*_+$ and $\|\bfz_t-\bfz^*\|=-t$. Thus, $\bfz_t\in \F \cap N(\bfz^*, \delta)$ for any $t\in (-\delta, 0)$. Applying \eqref{local-opt} yields
 $$g(\bfz_t)=t(\nabla_{\bfy} f(\bfz^*))_i+\lambda_{\max}t^2 \geq0,$$
which implies $(\nabla_{\bfy} f(\bfz^*))_i+\lambda_{\max}t  \leq 0.$ Letting $t\rightarrow0$ derives $(\nabla_{\bfy} f(\bfz^*))_i \leq 0$. This delivers $(\nabla_{\bfy} f(\bfz^*))_{\overline\Gamma^*} \leq 0$.

$\bullet$ $\|\bfx^*\|_0<s$.  The same reasoning for case  $\|\bfy^*_+\|_0<k$ leads to $\nabla_{\bfx} f(\bfz^*)= 0$.

$\bullet$  $\|\bfx^*\|_0=s$.  The same reasoning for case  $\|\bfy^*_+\|_0<k$ yields $(\nabla_{\bfx} f(\bfz^*))_{T^*}= 0$.
 
{Now we conclude that  case  $\|\bfy^*_+\|_0<k$ cannot occur. If we have  $\|\bfy^*_+\|_0<k$, then $2(A\bfx^* +\bfy^*-\epsilon {\bf1})=\nabla_{\bfy} f(\bfz^*)= 0$, resulting in $\nabla_{\bfx} f(\bfz^*)=2A^\top(A\bfx^* +\bfy^*-\epsilon {\bf1})+ 2\eta\bfx^*=2\eta\bfx^*=2\eta(\bfx^*_{T^*}; 0)$. This together with $(\nabla_{\bfx} f(\bfz^*))_{T^*}=0$ yields $\bfx^*_{T^*}=0$ and hence $\bfx^*=0$. Again by $A\bfx^* +\bfy^*-\epsilon {\bf1}= 0$, it follows $\bfy^*=\epsilon {\bf1}$ immediately. Apparently, this leads to a contradiction $k>\|\bfy^*_+\|_0=m\geq k$. Overall, we show  (\ref{KKT-point}).
}

{\bf b) Sufficiency.} Conversely, let $ \bfz^* $ satisfy (\ref{KKT-point}).  We  consider the following  three cases. 


$\bullet$ $\|\bfy^*_+ \|_0=k$. 
Consider a local region $N(\bfz^*, \delta_1)$ with $\delta_1:= \min \{ y_i^*: y_i^*> 0\}.$  {Thus, for any $ \bfz\in \F \cap N(\bfz^*, \delta_1)$, we have $y_j>0$ if $y_j^*>0$ and}
\begin{eqnarray}\label{yTyT}
y_{i}\leq0,~~\forall ~i\in \overline \Gamma^*.
\end{eqnarray}
In fact, if there exists an $i\in\overline \Gamma^*$ satisfying $y_{i}>0$, then  $\|\bfy_+\|_0\geq \|\bfy^*_+\|_0 +1=k+1$, which contradicts with $\bfy\in K$.  Direct calculations yield
\begin{eqnarray*}
&&\langle \nabla_{\bfy} f(\bfz^*),\bfy-\bfy^*\rangle\\
 &\overset{\eqref{KKT-point}}{=}&\langle 0,(\bfy-\bfy^*)_{\Gamma^*}\rangle+\langle ( \nabla_{\bfy} f(\bfz^*))_{\overline \Gamma^*},\bfy_{\overline \Gamma^*}\rangle\overset{\eqref{KKT-point}
 \eqref{yTyT}}{\geq}0.
\end{eqnarray*}

$\bullet$  $\|\bfx^* \|_0<s$.  This yields $\langle \nabla_{\bfx} f(\bfz^*),\bfx-\bfx^*\rangle=0$ due to $\nabla_{\bfx} f(\bfz^*)=0$ by \eqref{KKT-point}.

$\bullet$  $\|\bfx^* \|_0=s$.  Consider a local region $N(\bfz^*, \delta_2)$ with $\delta_2:=
\min \{ |x_i^*|: x_i^*\neq 0\}. $ For any $ \bfz:=(\bfx; \bfy)\in \F \cap N(\bfz^*, \delta_2)$,
$x_j \neq 0$ if $x_j^* \neq 0$, which indicates $T^*\subseteq\supp(\bfx)$. This together with $\bfx\in S$ that $\|\bfx\|_0\leq s=|T^*|$ suffices to
\begin{eqnarray}\label{xTxT}\supp(\bfx)= T^*.\end{eqnarray}
Now, one can verify that
\begin{eqnarray*}
&&\langle \nabla_{\bfx} f(\bfz^*),\bfx-\bfx^*\rangle\\
 &\overset{(\ref{KKT-point},\ref{xTxT})}{=}&\langle 0,(\bfx-\bfx^*)_{ T^*}\rangle+\langle ( \nabla_{\bfx} f(\bfz^*))_{\overline T^*},0\rangle=0.
\end{eqnarray*}
Using the above facts, for any $\bfz \in\F \cap N(\bfz^*, \min\{\delta_1,\delta_2\})$, it follows from \eqref{ell-0-1-obj-convex-smooth} that
\begin{eqnarray*}
&&f(\bfz)-f( \bfz^*)-  {\|\bfz - \bfz^*\|^2_{H}}\\
&=&\langle \nabla_{\bfx} f(\bfz^*),\bfx-\bfx^*\rangle + \langle \nabla_{\bfy} f(\bfz^*),\bfy-\bfy^*\rangle {\geq} 0.
\end{eqnarray*}
Therefore, $\bfz^*$ is the unique global minimizer of $\min  \{ f(\bfz):  \bfz\in \F \cap N(\bfz^*, \min\{\delta_1,\delta_2\})\}
$,  {which means it is} the unique local minimizer of problem  (\ref{ell-0-1}).
\end{IEEEproof}

\subsection{Proof of \Cref{pro-eta} }
\begin{IEEEproof}
 A $\tau $-stationary point satisfies \eqref{eta-point-z} which is equivalent to
  \begin{eqnarray}\label{eta-point}\begin{cases}
\bfx^* \in& \P_S\left( \bfx^* - \tau   \nabla_{\bfx} f(\bfz^*)\right),\\
\bfy^* \in& \P_K\left( \bfy^* - \tau   \nabla_{\bfy} f(\bfz^*)\right).
\end{cases}
  \end{eqnarray}
  Therefore, we show the equivalence between \eqref{eta-point-1} and \eqref{eta-point}. For $\bfx^*$ part, this can be guaranteed by \cite[Lemma 2.2]{beck2013sparsity}. For $\bfy^*$ part, projection  \eqref{pky}  enables to show that \eqref{eta-point-1} $\Rightarrow$ \eqref{eta-point}. So we only prove \eqref{eta-point} $\Rightarrow$ \eqref{eta-point-1}. Let $\bfl^*:=\nabla_{\bfy} f(\bfz^*)$. It follows from \eqref{pky} that
  \begin{eqnarray*}
  \bfy^*&\in&\P_K({\bfy}^*-\tau \bfl^*)\\
  &=&\left\{\left[
\begin{array}{c}
\bfy_{\Gamma}^*-\tau \bfl^*_{\Gamma}\\
0
\end{array}
\right]:~\Gamma\in\T(\bfy^*-\tau \bfl^*;k)\right\}.
  \end{eqnarray*}
This derives $\|\bfy^*_+\|_0\leq k$, and for any $\Gamma\in\T(\bfy^*-\tau\bfl^*;k)$,
\begin{eqnarray}\label{pro-z-T}\bfy^*_{\overline{\Gamma}}=0,~~~\bfl^*_\Gamma=0,~~~\bfy^*-\tau\bfl^*=\left[
\bfy^*_\Gamma;~-\tau\bfl^*_{\overline{\Gamma}}
\right], \end{eqnarray}
which together with the definition of $\T(\bfy^*-\tau\bfl^*;k)$ in \eqref{T-z} gives rise to
\begin{eqnarray*}
\Gamma =\Gamma_k   \cup \Gamma_-=\supp(\bfy^*),~~\overline \Gamma =(\Gamma_+\setminus\Gamma_{k})\cup \Gamma_0,   \end{eqnarray*}
where $\Gamma_+, \Gamma_-$ and $\Gamma_0$ are defined as \eqref{notation-z} in which $\bfy$ is replaced by $\bfy^*-\tau\bfl^*$. On index set $\overline{\Gamma}$, all elements $y^*_i-\tau\lambda^*_i=-\tau \lambda^*_{i}\geq 0$, namely, $\lambda^*_{i}\leq 0, i\in\overline{\Gamma}$.

{Now we claim that case $\|\bfy^*_+\|_0<k$ cannot happen. In fact, if $\|\bfy^*_+\|_0<k$}  and suppose there is an $i\in \overline{\Gamma}$ such that $\lambda^*_i<0$, then $\bfy^*-\tau\bfl^*$ has at least $\|\bfy^*_+\|_0+1\leq k$ positive entries and thus $$\|\bfy^*_+\|_0=\|(\P_K({\bfy}^*-\tau \bfl^*))_+\|_0 \geq \|\bfy^*_+\|_0+1.$$ This is a contradiction. So,  $\bfl^*_{ \overline{\Gamma}}=0$, { leading to $\bfl^*=\nabla_{\bfy} f(\bfz^*)=0$ by \eqref{pro-z-T}. Then, same reasoning to show the necessity of \Cref{nec-suff-opt-con-KKT} allows us to claim that case $\|\bfy^*_+\|_0<k$ cannot happen.  Therefore, we have $\|\bfy^*_+\|_0=k$. Then \eqref{eta-point-1} is satisfied  for any $j\in\supp(\bfy^*)=  \Gamma$ due to  $\bfl^*_{ \Gamma}=0$.}   For $j\notin\supp(\bfy^*)$, namely, $j\in \overline  \Gamma$,   the definition of $\Gamma_k$ in \eqref{T-z} yields
\begin{eqnarray*}
  0\leq y^*_j-\tau  \lambda^*_j \leq  y^*_i -\tau  \lambda^*_i , ~~\forall~i\in  \Gamma_k,
\end{eqnarray*}
which together with  $\Gamma_k\subseteq \Gamma$ and \eqref{pro-z-T} results in
\begin{eqnarray*}
 0\leq-\tau  \lambda^*_j\leq y^*_i, ~~ \forall~i\in \Gamma_k.
\end{eqnarray*}
Hence,  $-\|\bfy^*_+\|_{[k]}=-{\min}_{i\in\Gamma_k}y^*_i\leq\tau \lambda^*_j\leq 0 , \forall~j\in \overline  \Gamma (j\notin\supp(\bfy^*))$, showing \eqref{eta-point-1}.
\end{IEEEproof}

\subsection{Proof of \Cref{nec-suff-opt-con-sta}}
%

\begin{IEEEproof} Let $\bfz^*$ be a global minimizer. If it is not a $\tau$-stationary point with $0<\tau \leq 1/(2\lambda_{\max})$, then we have condition  $
\overline\bfz (\neq\bfz^*) \in \P_\F\left( \bfz^* - \tau \nabla  f(\bfz^*)\right).$ {Thus,}
  \begin{eqnarray*}
\|\overline\bfz-( \bfz^* - \tau   \nabla  f(\bfz^*))\|^2 < \|\bfz^*-( \bfz^* - \tau   \nabla f(\bfz^*))\|^2,
  \end{eqnarray*}
  which suffices to $$ 2\tau\langle   \nabla f(\bfz^*), \overline\bfz- \bfz^*\rangle < -  \|\overline\bfz-  \bfz^*  \|^2.$$
   {T}ogether with  \eqref{ell-0-1-obj-convex-smooth} and $0<\tau \leq 1/(2\lambda_{\max})$ derives
\begin{eqnarray*}
f( \overline\bfz) -  f( \bfz^*)&\leq&
  \langle \nabla  f(\bfz^*), \overline\bfz -\bfz^*\rangle
+ \lambda_{\max} \| \overline\bfz -\bfz^*\|^2\\
&<&( \lambda_{\max}-1/(2\tau))\| \overline\bfz-\bfz^*\|^2
 \leq 0.
\end{eqnarray*}
It contradicts the global optimality of $\bfz^*$. Therefore, $\bfz^*$ is a $\tau$-stationary point  with $0<\tau\leq 1/(2\lambda_{\max})$.

Conversely, let $\bfz^*$ be a $\tau$-stationary point  with $\tau \geq 1/(2\lambda_{\min})$. The definition  of $\P_\F$ and \eqref{eta-point-z} imply
   \begin{eqnarray*} \| \bfz^*-( \bfz^* - \tau  \nabla  f(\bfz^*))\|^2 \leq \|\bfz -( \bfz^* - \tau  \nabla  f(\bfz^*))\|^2,
  \end{eqnarray*}
for any $\bfz\in \F$, delivering $2\tau\langle  \nabla  f(\bfz^*),  \bfz- \bfz^*\rangle \geq -  \| \bfz-  \bfz^*  \|^2$.
This and \eqref{ell-0-1-obj-convex-smooth} yield
  \begin{eqnarray*}
f(\bfz ) - f( \bfz^*)
 &\geq &\langle \nabla f(\bfz^*),\bfz-\bfz^*\rangle  + \lambda_{\min}  \|\bfz-\bfz^*\|^2\\
 &\geq& (\lambda_{\min}  -1/
 (2\tau)) \|\bfz-\bfz^*\|^2.
\end{eqnarray*}
Since $\tau \geq 1/(2\lambda_{\min})$, the above relation shows the global optimality of $\bfz^*$ to \eqref{ell-0-1}. \end{IEEEproof}

\subsection{Proof of \Cref{check-global} }
\begin{IEEEproof}{Since $\bfz^*$ is a local minimizer of   problem (\ref{ell-0-1}),   it satisfies \eqref{KKT-point}. Then it follows 
\begin{eqnarray}\label{KKT-point-y1} 
    (\nabla_{\bfy} f(\bfz^*))_{\overline\Gamma^*} &=&2(A_{\overline\Gamma^* T^*}\bfx^*_{T^*}-\epsilon{\bf1})\leq 0, \\\label{KKT-point-y2} 
 (\nabla_{\bfy} f(\bfz^*))_{ \Gamma^*} &=&2(A_{ \Gamma^* T^*}\bfx^*_{T^*}+ \bfy^*_{ \Gamma^*} -\epsilon{\bf1})= 0, 
   \end{eqnarray}
   and thus 
\begin{eqnarray}   
\label{KKT-point-x2} 
(\nabla_{\bfx} f(\bfz^*))_{\overline T^*}   
&\overset{\eqref{KKT-point-y2}}{=}&   2A_{\overline\Gamma^* \overline T^*}^\top(A_{\overline\Gamma^* T^*}\bfx^*_{T^*}-\epsilon{\bf1}),\\
(\nabla_{\bfx} f(\bfz^*))_{T^*} 
&\overset{\eqref{KKT-point-y2}}{=}&  2A_{\overline\Gamma^* T^*}^\top(A_{\overline\Gamma^* T^*}\bfx^*_{T^*}-\epsilon{\bf1})+2\eta\bfx^*_{T^*}\nonumber\\
\label{tau-point-x2} &\overset{\eqref{KKT-point}}{=}&0.  
   \end{eqnarray}
   Using  last equation \eqref{tau-point-x2}  results in
   \begin{eqnarray*} 
(A_{\overline\Gamma^* T^*}A_{\overline\Gamma^* T^*}^\top+\eta I)(A_{\overline\Gamma^* T^*}\bfx^*_{T^*}-\epsilon{\bf1})=-\eta\epsilon{\bf1}. 
   \end{eqnarray*}
   This derives that
     \begin{eqnarray}\label{KKT-point-xi}   
    (\nabla_{\bfy} f(\bfz^*))_{\overline\Gamma^*} &\overset{\eqref{KKT-point-y1}}{=}&2(A_{\overline\Gamma^* T^*}\bfx^*_{T^*}-\epsilon{\bf1})\nonumber\\
    &=&-2\eta\epsilon(A_{\overline\Gamma^* T^*}A_{\overline\Gamma^* T^*}^\top+\eta I)^{-1}{\bf1}\nonumber\\
    &=&-\bxi^*.
   \end{eqnarray} 
  Then it is easy to see that $\bxi^*\geq 0$ by \eqref{KKT-point-y1} and 
       \begin{eqnarray*}  
 \forall~i\in\overline\Gamma^*:~~   -\tau_*  (\nabla_{\bfy} f(\bfz^*))_i &\overset{\eqref{KKT-point-xi}}{=}& \tau_*\xi^*_i \overset{\eqref{eta-point-global-tau*}}{\leq} \|\bfy^*_+\|_{[k]}.\end{eqnarray*} 
 If  $\|\bfx^*\|_0=s$ then $\|\bfx^*\|_{[s]}>0$ and thus for $\forall~i\in\overline T^*$\begin{eqnarray*} 
 \tau_* | (\nabla_{\bfx} f(\bfz^*))_i|&\overset{\eqref{KKT-point-x2}}{=}& \tau_*  |(A_{\overline\Gamma^* \overline T^*}^\top \bxi^*)_i|  \overset{\eqref{eta-point-global-tau*}}{\leq}  \|\bfx^*\|_{[s]}.
   \end{eqnarray*} 
If $\tau_* \geq 1/(2\lambda_{\min})$, then $\bfz^*$ is also a global minimizer by \Cref{nec-suff-opt-con-sta}.    The whole proof is finished.}  
\end{IEEEproof}

\section{Proofs of all theorems in Section \ref{sec:rec-per}}
\subsection{Proof of \Cref{feasible-recovery} }
\begin{IEEEproof}{
Conditions $\|\widehat{\bfx}\|=1$ and $ \|\widehat{\bfx}\|_0\leq s$ are clearly true. We only need to show that $\|(-A\widehat{\bfx})_+\|_0\leq k$, which is equivalent to prove that  $\|(-A \bfx^*)_+\|_0\leq k$.  For notational simplicity, denote $B:=A_{\overline\Gamma^* T^*}$ and $D:=(B B^\top+\eta I) $.  By noting that $A \bfx^* = (A_{\Gamma^*T^*}\bfx^*_{T^*}; B\bfx^*_{T^*})$, we next show that $$ \|(-A_{\Gamma^*T^*}\bfx^*_{T^*})_+\|_0\leq k,\quad B\bfx^*_{T^*}\geq0,$$ as it can claim $\|(-A \bfx^*)_+\|_0\leq k$ immediately.}

{Since $\bfz^*$ is a local minimizer of problem (\ref{ell-0-1}), it satisfies \eqref{eta-point-1} with $\tau=\tau_*$, thereby leading to (\ref{KKT-point-y1}-\ref{tau-point-x2}). From condition \eqref{KKT-point-y2}, i.e.,
   $\bfy^*_{ \Gamma^*} = \epsilon{\bf1} - A_{ \Gamma^* T^*}\bfx^*_{T^*}$, it follows $\|(\epsilon{\bf1} - A_{ \Gamma^* T^*}\bfx^*_{T^*})_+\|_0=  k$ due to $\|(\bfy^*_{ \Gamma^*})_+\|_0= \|\bfy^*_+\|_0=k$, thereby yielding $$\|( - A_{ \Gamma^* T^*}\bfx^*_{T^*})_+\|_0 \leq k.$$
 Using equation \eqref{tau-point-x2}  results in
   \begin{eqnarray*} 
   \bfx^*_{T^*} &=& \epsilon  
(B^\top B+\eta I)^{-1} B^\top {\bf1}\\
 &=& (\epsilon/\eta) \left(I- B^\top D^{-1} B\right) B^\top {\bf1}\\
 &=& (\epsilon/\eta) B^\top \left(I- D^{-1} BB^\top \right)  {\bf1}\\
 &=& (\epsilon/\eta) B^\top \left(I- D^{-1} (BB^\top+\eta I-\eta I) \right)  {\bf1}\\
&=& \epsilon B  ^\top 
D^{-1} {\bf1},
   \end{eqnarray*}
   where the second equation is due to the Woodbury matrix identity, namely,
$\left(I+UCV\right)^{-1}=I- U\left(C^{-1}+VU\right)^{-1}V.$
Therefore, we obtain
\begin{eqnarray*} 
\epsilon{\bf1} -  B\bfx^*_{T^*}
&=& \epsilon{\bf1} -  \epsilon BB  ^\top 
D^{-1} {\bf1}  \\
&=&\epsilon{\bf1} -  \epsilon (BB  ^\top+\eta I-\eta I) 
D^{-1} {\bf1}  \\
&=& \epsilon  \eta  
D^{-1} {\bf1}  \\ 
&=& \epsilon\eta   U  
{\rm Diag}\left( ({\lambda_1+\eta})^{-1} \cdots  ({\lambda_r+\eta})^{-1}\right)
  U^\top 
 {\bf1},
   \end{eqnarray*}
   where the last equation used the eigenvalue decomposition of $B B^\top =U{\rm Diag}(\lambda_1 \cdots \lambda_r)U^\top$ with $\lambda_1\geq\cdots\geq\lambda_r$ and $r=|\overline\Gamma^*|$. We note that $\lambda_r\geq c_A\geq0$.  If $\lambda_r=0$, then $\eta=0$ by \eqref{eta-tuning} and thus $\epsilon{\bf1} -  B\bfx^*_{T^*}=0$. This shows $B\bfx^*_{T^*} = \epsilon{\bf1}>0$. Now we focus on case $\lambda_r>0$. It is easy to see that $r=m-|\Gamma^*|\leq m-\|\bfy^*_+\|_0 = m-k$. Using these facts allows us to derive that  \begin{eqnarray*} \label{Bx-1-lbd}
0\overset{\eqref{KKT-point-y1}}{\leq} \epsilon{\bf1} -   B\bfx^*_{T^*}
 & {\leq}& \|\epsilon{\bf1} -   B\bfx^*_{T^*}\|_{\infty}{\bf1}   \nonumber\\
 & {\leq}& \|\epsilon{\bf1} -   B\bfx^*_{T^*}\| {\bf1}   \nonumber\\
&\leq&  \epsilon  {\eta }/({\lambda_r+\eta})  \| U^\top 
 {\bf1}  \|{\bf1}\nonumber\\
 &=&  \epsilon  {\eta }\sqrt{r}/({\lambda_r+\eta})  {\bf1}\nonumber\\
  &\leq&  \epsilon \eta \sqrt{m-k}/( {\lambda_r+\eta} ) {\bf1}\nonumber\\
&\overset{\eqref{eta-tuning}}{\leq}&  \epsilon{\bf1},
   \end{eqnarray*}
resulting in $ B\bfx^*_{T^*} \geq 0$ as well. } 
\end{IEEEproof} 
   
\subsection{Proof of Theorem \ref{oracle-property}}
 \begin{IEEEproof}{ 
 It follows from \cite[Lemma 4]{jacques2013robust} that
 \begin{eqnarray} \label{upper-bound-1}
\frac{1}{m}\| \sgn(\Phi\bfx^{\tr}+\var)-\sgn(\Phi\bfx^{\tr})\|_0  \leq \frac{\varrho}{2}+\delta
\end{eqnarray}
 with a probability higher than $1-e^{-2m\delta^2}$. Since $\Phi$ is generated as Lemma \ref{lemma-bse} with $m\geq m(2s,\delta,\theta)$ for a fixed $\theta\in[0,1]$ and $\delta\in(0,1)$, mapping $\sgn(\cdot)$  is then a  B$\delta$SE of order $2s$ with probability exceeding $1-\theta$ by  Lemma \ref{lemma-bse}. Hence,
\begin{eqnarray}\label{upper-bound-2}
d(\widehat \bfx, \bfx^{\tr}) 
     {\leq}   \frac{1}{m}\|\sgn(\Phi\widehat \bfx)- \sgn(\Phi\bfx^{\tr})\|_0 + \delta.\end{eqnarray}
We have shown  $\|(-A\widehat{\bfx})_+\|_0\leq k$ in \eqref{hat-x}, which together with $A={\rm Diag}(\bfc)\Phi$ implies $\|\sgn(\Phi\widehat \bfx)-\bfc\|_0\leq k$. Moreover, it is easy to see that $\|\bfc- \sgn(\Phi\bfx^{\tr}+\var)\|_0 \leq k $. Overall, using these facts  allows us to derive that
 {\allowdisplaybreaks  \begin{eqnarray} \label{fact-cos}
   && d(\widehat \bfx, \bfx^{\tr})\nonumber\\
    &\overset{\eqref{upper-bound-2}}{\leq}&  \frac{1}{m}\|\sgn(\Phi\widehat \bfx)- \sgn(\Phi\bfx^{\tr})\|_0 + \delta\nonumber\\
  &\leq&  \frac{1}{m}\|\sgn(\Phi\widehat \bfx)- \sgn(\Phi\bfx^{\tr}+\var)\|_0 \nonumber\\
  &+& \frac{1}{m}\| \sgn(\Phi\bfx^{\tr}+\var)- \sgn(\Phi\bfx^{\tr})\|_0  + \delta \nonumber\\
    &\leq&   \frac{1}{m}\|\sgn(\Phi\widehat \bfx)- \bfc\|_0 + \frac{1}{m}\|\bfc- \sgn(\Phi\bfx^{\tr}+\var)\|_0 \nonumber\\
  &+& \frac{1}{m}\| \sgn(\Phi\bfx^{\tr}+\var)- \sgn(\Phi\bfx^{\tr})\|_0  + \delta\nonumber\\
   &\overset{\eqref{upper-bound-1}}{\leq}& \frac{k}{m}  +\frac{k}{m} +\frac{\varrho}{2}+  2\delta =2\phi.
\end{eqnarray}}
Finally, it is easy to see that,  for any $\|\bfu\|=\|\bfv\|=1$,
\begin{eqnarray*}
d( \bfu,\bfv) =\frac{1}{\pi}{\rm arccos}\langle\bfu,\bfv\rangle=\frac{1}{\pi}{\rm arccos} \frac{2-\|\bfu-\bfv\|^2}{2}
\end{eqnarray*}
due to $0\leq\phi\leq1/2$, which suffices to
 \begin{eqnarray*}
d( \bfu,\bfv) \leq 2\phi ~ \Longrightarrow ~ \|\bfu-\bfv\|^2 &\leq& 2-2{\rm cos}(2\phi\pi)\\
&=& 2-2(1-2{\rm sin}^2(\phi\pi))\\
&=&4{\rm sin}^2(\phi\pi).
\end{eqnarray*}
This together with \eqref{fact-cos} displays the desired result. }
\end{IEEEproof}

 \section{Proofs of all theorems in Section \ref{sec:Newton}}
\subsection{Proof of \Cref{alpha-well-defined}}
 \begin{IEEEproof}It follows from $\bfz^\ell(\tau)  \in \P_{\F}(\bfz^\ell -\tau \nabla  f(\bfz^\ell))$ that
$$\| \bfz^\ell(\tau)-( \bfz^\ell - \tau  \nabla  f(\bfz^\ell))\|^2 \leq \|\bfz^\ell -( \bfz^\ell - \tau  \nabla  f(\bfz^\ell))\|^2,$$
 which results in
   \begin{eqnarray} \label{eta-point-not-3}
 2\tau\langle   \nabla  f(\bfz^\ell),  \bfz^\ell(\tau)-  \bfz^\ell\rangle \leq -  \| \bfz^\ell(\tau)-  \bfz^\ell\|^2.
  \end{eqnarray}
 {Combining with} \eqref{ell-0-1-obj-convex-smooth} leads to
\begin{eqnarray*}
&& f(\bfz^\ell(\tau))\\
  &\leq&    f(\bfz^\ell)+
  \langle  \nabla  f(\bfz^\ell),\bfz^\ell(\tau)-\bfz^\ell\rangle
+ \lambda_{\max}\|\bfz^\ell(\tau)-\bfz^\ell\|^2\\
& {\leq}&   f(\bfz^\ell)-\Big( {1}/{(2\tau)} -\lambda_{\max}\Big)\|\bfz^\ell(\tau)-\bfz^\ell\|^2\\
&\leq&  f(\bfz^\ell) -  \rho\|\bfz^\ell(\tau)-\bfz^\ell\|^2,
\end{eqnarray*}
where the last inequality is from $0<\tau\leq  1/(2\rho+2\lambda_{\max})$.  {Invoking} the Armijo-type step size rule,  {one has} $\tau_\ell\geq \beta/(2\rho+2\lambda_{\max})$,  {which by  $\tau_\ell\leq 1$ proves the desired assertion.}
\end{IEEEproof}

\subsection{Proof of \Cref{zk1-zk-0}}
\begin{IEEEproof}
i) By \Cref{alpha-well-defined} and $\bfu^\ell=\bfz^{\ell}(\tau_\ell)$, we have
\begin{eqnarray}\label{u-z-ell}
		 f(\bfu^{\ell})&\leq&  f(\bfz^\ell) -  \rho \|\bfu^\ell-\bfz^\ell\|^2.
		\end{eqnarray}
By the framework of \Cref{Alg-NL01}, if $\bfz^{\ell+1}=\bfu^\ell$, then the above condition implies,
\begin{eqnarray}
\label{fact-u-z}		 f(\bfz^{\ell+1})&\leq&
		  f(\bfz^\ell) -  \rho \|\bfu^{\ell}-\bfz^\ell\|^2\\
		  &=&
		  f(\bfz^\ell) - \rho \|\bfz^{\ell+1}-\bfz^\ell\|^2.\nonumber
		\end{eqnarray}
If $\bfz^{\ell+1}=\bfv^\ell$, then  we obtain
\begin{eqnarray}
		f(\bfz^{\ell+1})&=& f(\bfv^{\ell})\nonumber\\
		&\leq&f(\bfu^\ell) - \rho\|\bfz^{\ell+1} -\bfu^\ell\|^2\nonumber\\
\label{fact-u-z-1}			&\leq& f(\bfz^\ell) - \rho\|\bfu^\ell-\bfz^\ell\|^2- \rho\|\bfz^{\ell+1}-\bfu^\ell\|^2~~~~\\
		&\leq& f(\bfz^\ell) - (\rho/2) \|\bfz^{\ell+1}-\bfz^\ell\|^2, \nonumber
		\end{eqnarray}
		where the second and last inequalities used \eqref{u-z-ell} and a fact $\|{\bf a}+\bfb\|^2\leq 2\|{\bf a}\|^2+2\|\bfb\|^2$  {for all vectors ${\bf a}$ and $\bfb$}. Both cases lead to	 \begin{eqnarray}\label{decrese-pro}
		 f(\bfz^{\ell+1})&\leq& f(\bfz^\ell) - (\rho/2)\|\bfz^{\ell+1}-\bfz^\ell\|^2, \\
		  f(\bfz^{\ell+1})&\leq& f(\bfz^\ell) -  \rho \|\bfu^{\ell}-\bfz^\ell\|^2. \nonumber
		\end{eqnarray}
Therefore, $\{f(\bfz^{\ell})\}$ is a non-increasing sequence, and thus
$$ \max\{ \| A\bfx^{\ell}-\epsilon {\bf1}+\bfy^{\ell}\|^2,  \eta \|\bfx^{\ell}\|^2\} \leq f(\bfz^{\ell})\leq f(\bfz^{0}),$$
which  {indicates the boundedness of} $\{ \bfx^{\ell}\}$ and  $\{ \bfy^{\ell}\}$, and so that of  $\{ \bfz^{\ell}\}$. The non-increasing property in \eqref{decrese-pro}
 and $f \geq 0$ also give rise to
		\begin{eqnarray*} &&\sum_{\ell\geq 0} \max\{(\rho/2)\|\bfz^{\ell+1}-\bfz^\ell\|^2, \rho\|\bfu^{\ell}-\bfz^\ell\|^2\}\\
		&\leq& \sum_{\ell\geq 0} \left[ f(\bfz^\ell) - f(\bfz^{\ell+1})\right]\\
		&=& f(\bfz^0) - \lim_{\ell\rightarrow\infty} f(\bfz^{\ell+1})\leq f(\bfz^0).	
		\end{eqnarray*}
The above condition suffices to $\lim_{\ell\rightarrow\infty}\|\bfz^{\ell+1}-\bfz^\ell\|=\lim_{\ell\rightarrow\infty}\|\bfu^{\ell}-\bfz^\ell\|=0.$

ii) Let $\bfz^*$ be any accumulating point of $\{\bfz^\ell\}$.  {Then there exists a} subset $J$ of $\{0,1,2,\cdots\}$ such that $\lim_{\ell(\in J)\rightarrow\infty} \bfz^\ell = \bfz^*.$  {This further implies $\lim_{\ell(\in J)\rightarrow\infty} \bfu^\ell = \bfz^*$ by applying a). 
In addition, as stated in \Cref{alpha-well-defined}, we have $\{\tau_\ell\}\subseteq [ \underline{\tau}, 1]$,
which indicates that one can find a subsequence $L$ of $J$ and a scalar $\tau_*\in [ \underline{\tau}, 1]$ such that $\{\tau_\ell: \ell\in L\}\rightarrow \tau_*$. 
To summarize, we have}
\begin{eqnarray}
		\label{a-a-under} 
		\bfz^\ell\rightarrow \bfz^*,~~ \bfu^\ell \rightarrow \bfz^*,~~\tau_\ell \rightarrow\tau_* \in [\underline{\tau},1],~~ \ell(\in L)\rightarrow\infty.\end{eqnarray}
Let $\overline\bfz^\ell :=\bfz^\ell -\tau_\ell \nabla  f(\bfz^\ell)$.  \Cref{Alg-NL01} implies
\begin{eqnarray}\label{u-l-P}
\bfu^{\ell}  \in \P_\F( \overline\bfz^\ell ),~~ \lim_{\ell(\in L)\rightarrow\infty} \overline\bfz^\ell =\bfz^* -\tau_*  \nabla f(\bfz^*)=:\overline\bfz^*.\end{eqnarray}
The first condition means $\bfu^{\ell}\in \F$ for any $\ell\geq1$. Note that $\F$ is closed and  $\bfz^*$ is the accumulating point of $\{\bfu^{\ell}\}$ by \eqref{a-a-under}. Therefore, $\bfz^*\in\F$, which results in
\begin{eqnarray}
 \label{z-*-g}
\min_{\bfz\in\F}\|\bfz -\overline\bfz^*\| \leq \|\bfz^*-\overline{\bfz}^*\|.\end{eqnarray}
If the strict inequality holds in the above condition, then there is an $ \varepsilon_0>0$ such that
\begin{eqnarray*}
\|\bfz^*-\overline\bfz^*\|-\varepsilon_0 &=&  \min_{\bfz\in\F}\|\bfz -\overline\bfz^*\| \\
&\geq&\min_{\bfz\in\F} (\|\bfz - \overline\bfz^\ell \|-\|\overline\bfz^\ell- \overline\bfz^*\|)\\
&=&\|\bfu^{\ell} - \overline\bfz^\ell \|-\|\overline\bfz^\ell -\overline\bfz^*\|
\end{eqnarray*}
where the last equality is from \eqref{u-l-P}. Taking the limit of both sides of the above condition along $\ell(\in L)\rightarrow\infty$ yields $\|\bfz^*- \overline\bfz^*\|-\varepsilon_0 \geq \|\bfz^* - \overline\bfz^*\|$ by \eqref{a-a-under} and \eqref{u-l-P}, a contradiction with  $ \varepsilon_0>0$. Therefore, we must have the equality holds in \eqref{z-*-g}, showing that
\begin{eqnarray*} \bfz^* \in  \P_\F(\overline\bfz^*)= \P_\F\left( \bfz^* - \tau_*   \nabla f(\bfz^*)\right). \end{eqnarray*}
The above relation means the conditions in  \eqref{eta-point-1} hold  for $\tau=\tau_*$, then these conditions must hold for any $0<\tau\leq \underline{\tau}$ due to $\underline{\tau}\leq \tau_*$ from \eqref{a-a-under}, namely,
\begin{eqnarray*} \bfz^* \in   \P_\F\left( \bfz^* - \tau \nabla f(\bfz^*)\right), \end{eqnarray*}
displaying that $\bfz^*$ is a $\tau$-stationary point of \eqref{ell-0-1}.
\end{IEEEproof}
\subsection{Proof of \Cref{global-convergence}}
\begin{IEEEproof}As shown in \Cref{zk1-zk-0}, one can find a subsequence of $\{\bfz^\ell\}$ that converges to the $\tau$-stationary point $\bfz^*$ with $0<\tau\leq \underline{\tau}$ of \eqref{ell-0-1}. 
Recall that a $\tau $-stationary point $\bfz^*$ that satisfies \eqref{eta-point-1} also meets \eqref{KKT-point}, {which by \Cref{nec-suff-opt-con-KKT} indicates that $\bfz^*$ is the unique local minimizer.} In other words, $\bfz^*$ is an isolated local minimizer of (\ref{ell-0-1}).  Finally, it follows from $\bfz^*$ being isolated, \cite[Lemma 4.10]{more1983computing} and $\lim_{\ell\rightarrow\infty}\|\bfz^{\ell+1}-\bfz^\ell\|=0$ by \Cref{zk1-zk-0} that
the whole sequence converges to $\bfz^*$.  { It is also a  global minimizer due to Corollary \ref{check-global}.}
\end{IEEEproof}

\subsection{Proof of \Cref{convergence-rate} }

\begin{IEEEproof}
 {It follows Theorem \ref{global-convergence} and Lemma \ref{zk1-zk-0} that both $\bfz^\ell(\in\F)\rightarrow \bfz^* $ and $\bfu^\ell(\in\F)\rightarrow \bfz^*$ and $\bfz^*$ is a \ts. Thus $\|  \bfy^*_+\|_0 = k$. This implies  $\Gamma_+^\ell\equiv \widetilde \Gamma_+^\ell\equiv\Gamma_+^*$ for sufficiently large $\ell$. If $\|\bfx^* \|_0 = s$, then $\bfz^\ell(\in\F)\rightarrow \bfz^* $ indicates  $ T^{\ell}\equiv \widetilde  T^{\ell}\equiv T^*$  for sufficiently large $\ell$. If $\|\bfx^* \|_0 < s$, then by \eqref{eta-point-1} that $\nabla_{\bfx} f(\bfz^*)=0$. This together with $\bfu^\ell(\in\F)\rightarrow \bfz^* $ results in $\|\nabla_{\bfx} f(\bfu^\ell)\|<\varepsilon$. Therefore,  the framework of \Cref{Alg-NL01} allows us to assert that  $\bfz^{\ell+1}=\bfv^{\ell}$  for all sufficiently large $\ell$ due to $0<\rho\leq \lambda_{\min}$ and
\begin{eqnarray*}
  f(\bfv^{\ell})&\overset{\eqref{newton-step-v}}{\leq}&   f(\bfu^{\ell})-\lambda_{\min} \|\bfu^{\ell}-\bfv^{\ell}\|^2\\
  &\leq&   f(\bfu^{\ell})-\rho\|\bfu^{\ell}-\bfv^{\ell}\|^2.
\end{eqnarray*}
The above assertion means \eqref{Newton-descent-property} is always admitted for sufficiently large $\ell$. However, updating rule  \eqref{Newton-descent-property}  of  $\bfv^{\ell}$ indicates
\begin{eqnarray*}
T^{\ell+1}\subseteq T^{\ell},~~ \Gamma_+^{\ell+1} \subseteq \Gamma_+^\ell,~~ \Gamma_0^{\ell+1} \supseteq \Gamma_0^\ell,~~ \Gamma_-^{\ell+1}  {\subseteq} \Gamma_-^\ell.
\end{eqnarray*}
 Note that these sets have finite elements. Therefore, the sequences $\{T^{\ell}\}$,  $\{\Gamma_+^\ell\}$, $\{\Gamma_0^\ell\}$, $\{\Gamma_-^\ell\}$ converge. In other words,  there is a finite $\kappa\geq 1$ such that, for any $\ell\geq \kappa$,  
 \begin{eqnarray}
T^{\ell+1} = T^{\ell},~~ \Gamma_+^{\ell+1} = \Gamma_+^\ell,~~ \Gamma_0^{\ell+1} = \Gamma_0^\ell,~~ \Gamma_-^{\ell+1} = \Gamma_-^\ell,
\end{eqnarray}
which implies $\Omega(\bfz^{\ell+1})=\Omega(\bfz^{\ell})$.  This and \eqref{Newton-descent-property} yield
\begin{eqnarray*}
\bfz^{\ell+2}=\bfv^{\ell+1} &=&  {\rm argmin}~\{f\left(\bfz\right):\bfz\in \Omega(\bfz^{\ell+1})\}\\
	 &= & {\rm argmin}~\{f\left(\bfz\right):\bfz\in \Omega(\bfz^{\ell})\}\\
 &= & \bfv^{\ell}=\bfz^{\ell+1}.
\end{eqnarray*}
Overall, for any $\ell \geq \kappa$, we have $\bfz^{\ell} =\bfz^\kappa.$
Recall \Cref{global-convergence} that whole sequence $\{\bfz^{\ell}\}$ converges to $\bfz^*$, which suffices to $\bfz^*=\lim_{\ell\rightarrow\infty}\bfz^{\ell}=\bfz^\kappa.$ }
\end{IEEEproof}

 \subsection{Proof of \Cref{convergence-rate-true}}
\begin{IEEEproof}
{
Since limit $\bfz^{*}$ is a local minimizer by Theorem \ref{global-convergence}, we have \eqref{gap-x-xtrue} from Theorem \ref{oracle-property}, namely,  
   \begin{eqnarray*} 
\|c_*\bfx^*-\bfx^\tr\|  \leq 2{\rm sin}\left(\phi\pi\right). 
\end{eqnarray*}
Theorem \ref{convergence-rate}   ii) shows that \gpa\ will terminate at limit $\bfz^*$  within finite steps. That is, there is a finite $\kappa$ such that $ \bfz^\ell=\bfz^*$  for any $\ell\geq \kappa$, which derives 
\begin{eqnarray*}
\| c_*\bfx^{\ell}-\bfx^\tr\|= \| c_*\bfx^{*}-\bfx^\tr\|\leq 2{\rm sin}\left(\phi\pi\right).
\end{eqnarray*}
Hence, the whole proof is completed.}
\end{IEEEproof}
\section*{Acknowledgment}
The authors sincerely thank the associate editor and the five referees for their constructive comments, which have significantly improved the quality of the paper.

\ifCLASSOPTIONcaptionsoff
  \newpage
\fi



%


%
%
%
\bibliographystyle{IEEEtran}
\bibliography{references}

\begin{thebibliography}{10}
\providecommand{\url}[1]{#1}
\csname url@samestyle\endcsname
\providecommand{\newblock}{\relax}
\providecommand{\bibinfo}[2]{#2}
\providecommand{\BIBentrySTDinterwordspacing}{\spaceskip=0pt\relax}
\providecommand{\BIBentryALTinterwordstretchfactor}{4}
\providecommand{\BIBentryALTinterwordspacing}{\spaceskip=\fontdimen2\font plus
\BIBentryALTinterwordstretchfactor\fontdimen3\font minus
  \fontdimen4\font\relax}
\providecommand{\BIBforeignlanguage}[2]{{%
\expandafter\ifx\csname l@#1\endcsname\relax
\typeout{** WARNING: IEEEtran.bst: No hyphenation pattern has been}%
\typeout{** loaded for the language `#1'. Using the pattern for}%
\typeout{** the default language instead.}%
\else
\language=\csname l@#1\endcsname
\fi
#2}}
\providecommand{\BIBdecl}{\relax}
\BIBdecl

\bibitem{candes2005decoding}
E.~J. Candes and T.~Tao, ``Decoding by linear programming,'' \emph{IEEE Trans.
  Inf. Theory}, vol.~51, no.~12, pp. 4203--4215, 2005.

\bibitem{candes2006robust}
E.~J. Cand{\`e}s, J.~Romberg, and T.~Tao, ``Robust uncertainty principles:
  Exact signal reconstruction from highly incomplete frequency information,''
  \emph{IEEE Trans. Inf. Theory}, vol.~52, no.~2, pp. 489--509, 2006.

\bibitem{donoho2006compressed}
D.~L. Donoho, ``Compressed sensing,'' \emph{IEEE Trans. Inf. Theory}, vol.~52,
  no.~4, pp. 1289--1306, 2006.

\bibitem{BB08}
P.~T. Boufounos and R.~G. Baraniuk, ``1-bit compressive sensing,'' in
  \emph{42nd Annu. Conf. Inf. Sci. Syst.}\hskip 1em plus 0.5em minus
  0.4em\relax IEEE, 2008, pp. 16--21.

\bibitem{haboba2011architecture}
J.~Haboba, M.~Mangia, R.~Rovatti, and G.~Setti, ``An architecture for 1-bit
  localized compressive sensing with applications to {EEG},'' in \emph{IEEE
  Biomed. Circuits Syst. Conf.}\hskip 1em plus 0.5em minus 0.4em\relax IEEE,
  2011, pp. 137--140.

\bibitem{movahed2014iterative}
A.~Movahed and M.~C. Reed, ``Iterative detection for compressive sensing:
  {T}urbo {CS},'' in \emph{IEEE Int. Conf. Commun.}\hskip 1em plus 0.5em minus
  0.4em\relax IEEE, 2014, pp. 4518--4523.

\bibitem{tang2017low}
W.~Tang, W.~Xu, X.~Zhang, and J.~Lin, ``A low-cost channel feedback scheme in
  mmwave massive mimo system,'' in \emph{3rd IEEE Int. Conf. Intell. Comput.
  Commun.}\hskip 1em plus 0.5em minus 0.4em\relax IEEE, 2017, pp. 89--93.

\bibitem{zhou2017sparse}
Z.~Zhou, X.~Chen, D.~Guo, and M.~L. Honig, ``Sparse channel estimation for
  massive {MIMO} with 1-bit feedback per dimension,'' in \emph{IEEE Wirel.
  Commun. Netw. Conf.}\hskip 1em plus 0.5em minus 0.4em\relax IEEE, 2017, pp.
  1--6.

\bibitem{meng2009sparse}
J.~Meng, H.~Li, and Z.~Han, ``Sparse event detection in wireless sensor
  networks using compressive sensing,'' in \emph{43rd Annu. Conf. Inf. Sci.
  Syst.}\hskip 1em plus 0.5em minus 0.4em\relax IEEE, 2009, pp. 181--185.

\bibitem{feng2009multiple}
C.~Feng, S.~Valaee, and Z.~Tan, ``Multiple target localization using
  compressive sensing,'' in \emph{IEEE Glob. Commun. Conf.}\hskip 1em plus
  0.5em minus 0.4em\relax IEEE, 2009, pp. 1--6.

\bibitem{xiong20141}
J.~Xiong and Q.~Tang, ``1-bit compressive data gathering for wireless sensor
  networks,'' \emph{J. Sens.}, vol. 2014, 2014.

\bibitem{chen2015amplitude}
C.-H. Chen and J.-Y. Wu, ``Amplitude-aided 1-bit compressive sensing over noisy
  wireless sensor networks,'' \emph{IEEE Wirel. Commun. Lett.}, vol.~4, no.~5,
  pp. 473--476, 2015.

\bibitem{lee2012spectrum}
D.~Lee, T.~Sasaki, T.~Yamada, K.~Akabane, Y.~Yamaguchi, and K.~Uehara,
  ``Spectrum sensing for networked system using 1-bit compressed sensing with
  partial random circulant measurement matrices,'' in \emph{IEEE 75th Veh.
  Technol. Conf.}\hskip 1em plus 0.5em minus 0.4em\relax IEEE, 2012, pp. 1--5.

\bibitem{fu2014sub}
N.~Fu, L.~Yang, and J.~Zhang, ``Sub-nyquist 1 bit sampling system for sparse
  multiband signals,'' in \emph{22nd Eur. Signal Process. Conf.}\hskip 1em plus
  0.5em minus 0.4em\relax IEEE, 2014, pp. 736--740.

\bibitem{bourquard2012binary}
A.~Bourquard and M.~Unser, ``Binary compressed imaging,'' \emph{IEEE Trans.
  Image Process.}, vol.~22, no.~3, pp. 1042--1055, 2012.

\bibitem{dong2015map}
X.~Dong and Y.~Zhang, ``A {MAP} approach for 1-bit compressive sensing in
  synthetic aperture radar imaging,'' \emph{IEEE Geosci. Remote Sens. Lett.},
  vol.~12, no.~6, pp. 1237--1241, 2015.

\bibitem{marcos2016compressed}
D.~Marcos, T.~Lasser, A.~L{\'o}pez, and A.~Bourquard, ``Compressed imaging by
  sparse random convolution,'' \emph{Opt. Express}, vol.~24, no.~2, pp.
  1269--1290, 2016.

\bibitem{li2018survey}
Z.~Li, W.~Xu, X.~Zhang, and J.~Lin, ``A survey on one-bit compressed sensing:
  Theory and applications,'' \emph{Front. Comput. Sci.}, vol.~12, no.~2, pp.
  217--230, 2018.

\bibitem{qin2018sparse}
Z.~Qin, J.~Fan, Y.~Liu, Y.~Gao, and G.~Y. Li, ``Sparse representation for
  wireless communications: {A} compressive sensing approach,'' \emph{IEEE
  Signal Process. Mag.}, vol.~35, no.~3, pp. 40--58, 2018.

\bibitem{gopi2013one}
S.~Gopi, P.~Netrapalli, P.~Jain, and A.~Nori, ``One-bit compressed sensing:
  Provable support and vector recovery,'' in \emph{Int. Conf. Mach.
  Learn.}\hskip 1em plus 0.5em minus 0.4em\relax PMLR, 2013, pp. 154--162.

\bibitem{wanghuang2015}
H.~Wang, X.~Huang, Y.~Liu, H.~Sabine~Van, and W.~Qun, ``Binary reweighted
  l1-norm minimization for one-bit compressed sensing,'' in \emph{8th Int. Jt.
  Conf. Biomed. Eng. Syst. Technol.}, 2015.

\bibitem{flodin2019superset}
L.~Flodin, V.~Gandikota, and A.~Mazumdar, ``Superset technique for approximate
  recovery in one-bit compressed sensing,'' \emph{arXiv preprint
  arXiv:1910.13971}, 2019.

\bibitem{xiao2019one}
P.~Xiao, B.~Liao, and J.~Li, ``One-bit compressive sensing via schur-concave
  function minimization,'' \emph{IEEE Trans. Signal Process.}, vol.~67, no.~16,
  pp. 4139--4151, 2019.

\bibitem{khobahi2020model}
S.~Khobahi and M.~Soltanalian, ``Model-based deep learning for one-bit
  compressive sensing,'' \emph{IEEE Trans. Signal Process.}, vol.~68, pp.
  5292--5307, 2020.

\bibitem{laska2011trust}
J.~N. Laska, Z.~Wen, W.~Yin, and R.~G. Baraniuk, ``Trust, but verify: Fast and
  accurate signal recovery from 1-bit compressive measurements,'' \emph{IEEE
  Trans. Signal Process.}, vol.~59, no.~11, pp. 5289--5301, 2011.

\bibitem{needell2009cosamp}
D.~Needell and J.~A. Tropp, ``Co{S}a{MP}: Iterative signal recovery from
  incomplete and inaccurate samples,'' \emph{Appl. Comput. Harmon. Anal.},
  vol.~26, no.~3, pp. 301--321, 2009.

\bibitem{boufounos2009greedy}
P.~T. Boufounos, ``Greedy sparse signal reconstruction from sign
  measurements,'' in \emph{Conf. Rec. Asilomar Conf. Signals Syst.
  Comput.}\hskip 1em plus 0.5em minus 0.4em\relax IEEE, 2009, pp. 1305--1309.

\bibitem{jacques2013robust}
L.~Jacques, J.~N. Laska, P.~T. Boufounos, and R.~G. Baraniuk, ``Robust 1-bit
  compressive sensing via binary stable embeddings of sparse vectors,''
  \emph{IEEE Trans. Inf. Theory}, vol.~59, no.~4, pp. 2082--2102, 2013.

\bibitem{friedlander2020nbiht}
M.~P. Friedlander, H.~Jeong, Y.~Plan, and O.~Yilmaz, ``{NBIHT: A}n efficient
  algorithm for 1-bit compressed sensing with optimal error decay rate,''
  \emph{arXiv preprint arXiv:2012.12886}, 2020.

\bibitem{plan2012robust}
Y.~Plan and R.~Vershynin, ``Robust 1-bit compressed sensing and sparse logistic
  regression: A convex programming approach,'' \emph{IEEE Trans. Inf. Theory},
  vol.~59, no.~1, pp. 482--494, 2012.

\bibitem{zhang2014efficient}
L.~Zhang, J.~Yi, and R.~Jin, ``Efficient algorithms for robust one-bit
  compressive sensing,'' in \emph{Int. Conf. Mach. Learn.}\hskip 1em plus 0.5em
  minus 0.4em\relax PMLR, 2014, pp. 820--828.

\bibitem{cai2014soft}
X.~Cai, Z.~Zhang, H.~Zhang, and C.~Li, ``Soft consistency reconstruction: a
  robust 1-bit compressive sensing algorithm,'' in \emph{IEEE Int. Conf.
  Commun.}\hskip 1em plus 0.5em minus 0.4em\relax IEEE, 2014, pp. 4530--4535.

\bibitem{huang2018pinball}
X.~Huang, L.~Shi, M.~Yan, and J.~A. Suykens, ``Pinball loss minimization for
  one-bit compressive sensing: Convex models and algorithms,''
  \emph{Neurocomputing}, vol. 314, pp. 275--283, 2018.

\bibitem{rencker2019sparse}
L.~Rencker, F.~Bach, W.~Wang, and M.~D. Plumbley, ``Sparse recovery and
  dictionary learning from nonlinear compressive measurements,'' \emph{IEEE
  Trans. Signal Process.}, vol.~67, no.~21, pp. 5659--5670, 2019.

\bibitem{huang2018nonconvex}
X.~Huang and M.~Yan, ``Nonconvex penalties with analytical solutions for
  one-bit compressive sensing,'' \emph{Signal Process.}, vol. 144, pp.
  341--351, 2018.

\bibitem{xu2020feature}
W.~Xu, Y.~Tian, S.~Wang, and Y.~Cui, ``Feature selection and classification of
  noisy proteomics mass spectrometry data based on one-bit perturbed compressed
  sensing,'' \emph{Bioinformatics}, vol.~36, no.~16, pp. 4423--4431, 2020.

\bibitem{fu2014robust}
X.~Fu, F.-M. Han, and H.~Zou, ``Robust 1-bit compressive sensing against sign
  flips,'' in \emph{Glob. Commun. Conf.}\hskip 1em plus 0.5em minus 0.4em\relax
  IEEE, 2014, pp. 3121--3125.

\bibitem{li2014robust}
F.~Li, J.~Fang, H.~Li, and L.~Huang, ``Robust one-bit bayesian compressed
  sensing with sign-flip errors,'' \emph{IEEE Signal Process. Lett.}, vol.~22,
  no.~7, pp. 857--861, 2014.

\bibitem{shi2016methods}
H.-J.~M. Shi, M.~Case, X.~Gu, S.~Tu, and D.~Needell, ``Methods for quantized
  compressed sensing,'' in \emph{Inf. Theory Appl. Workshop}.\hskip 1em plus
  0.5em minus 0.4em\relax IEEE, 2016, pp. 1--9.

\bibitem{huang2018robust}
J.~Huang, Y.~Jiao, X.~Lu, and L.~Zhu, ``Robust decoding from 1-bit compressive
  sampling with ordinary and regularized least squares,'' \emph{SIAM J. Sci.
  Comput.}, vol.~40, no.~4, pp. A2062--A2086, 2018.

\bibitem{fan2021robust}
Q.~Fan, C.~Jia, J.~Liu, and Y.~Luo, ``Robust recovery in 1-bit compressive
  sensing via $l_q$-constrained least squares,'' \emph{Signal Process.}, vol.
  179, p. 107822, 2021.

\bibitem{yan2012robust}
M.~Yan, Y.~Yang, and S.~Osher, ``Robust 1-bit compressive sensing using
  adaptive outlier pursuit,'' \emph{IEEE Trans. Signal Process.}, vol.~60,
  no.~7, pp. 3868--3875, 2012.

\bibitem{movahed2012robust}
A.~Movahed, A.~Panahi, and G.~Durisi, ``A robust rfpi-based 1-bit compressive
  sensing reconstruction algorithm,'' in \emph{IEEE Inf. Theory
  Workshop}.\hskip 1em plus 0.5em minus 0.4em\relax IEEE, 2012, pp. 567--571.

\bibitem{movahed2014recovering}
A.~Movahed, A.~Panahi, and M.~C. Reed, ``Recovering signals with variable
  sparsity levels from the noisy 1-bit compressive measurements,'' in
  \emph{IEEE Int. Conf. Acoust. Speech Signal Process.}\hskip 1em plus 0.5em
  minus 0.4em\relax IEEE, 2014, pp. 6454--6458.

\bibitem{dai2016noisy}
D.-Q. Dai, L.~Shen, Y.~Xu, and N.~Zhang, ``Noisy 1-bit compressive sensing:
  models and algorithms,'' \emph{Appl. Comput. Harmon. Anal.}, vol.~40, no.~1,
  pp. 1--32, 2016.

\bibitem{armijo1966minimization}
L.~Armijo, ``Minimization of functions having lipschitz continuous first
  partial derivatives,'' \emph{Pac. J. Math.}, vol.~16, no.~1, pp. 1--3, 1966.

\bibitem{gianelli2016one}
C.~Gianelli, L.~Xu, J.~Li, and P.~Stoica, ``One-bit compressive sampling with
  time-varying thresholds: Maximum likelihood and the cram{\'e}r-rao bound,''
  in \emph{50th Asilomar Conf. Signals, Syst. Comput.}\hskip 1em plus 0.5em
  minus 0.4em\relax IEEE, 2016, pp. 399--403.

\bibitem{zhu2019computationally}
H.~Zhu, F.~Liu, and J.~Li, ``Computationally efficient sinusoidal parameter
  estimation from signed measurements: {ADMM} approaches,'' \emph{IEEE Signal
  Process. Lett.}, vol.~26, no.~12, pp. 1798--1802, 2019.

\bibitem{ren2019sinusoidal}
J.~Ren, T.~Zhang, J.~Li, and P.~Stoica, ``Sinusoidal parameter estimation from
  signed measurements via majorization--minimization based {RELAX},''
  \emph{IEEE Trans. Signal Process.}, vol.~67, no.~8, pp. 2173--2186, 2019.

\bibitem{beck2013sparsity}
A.~Beck and Y.~C. Eldar, ``Sparsity constrained nonlinear optimization:
  Optimality conditions and algorithms,'' \emph{SIAM J. Optim.}, vol.~23,
  no.~3, pp. 1480--1509, 2013.

\bibitem{more1983computing}
J.~J. Mor{\'e} and D.~C. Sorensen, ``Computing a trust region step,''
  \emph{SIAM J. Sci. Statist. Comput.}, vol.~4, no.~3, pp. 553--572, 1983.

\end{thebibliography}

%

%
%
%




\end{document}